\documentclass[a4paper]{article}
\usepackage{graphicx}
\usepackage{float}
\usepackage{amsmath}
\usepackage{amsthm}
\usepackage{amsfonts}
\usepackage{lmodern}
\usepackage{mathtools}
\usepackage{url}
\usepackage{hyperref}
\hypersetup{hidelinks}

\usepackage{xcolor}

\usepackage{geometry}
\usepackage{tikz}
\tikzset{
  myedge/.style={draw, black},
  mylabel/.style={midway, above, text=blue, fill=white, inner sep=1pt,
    font=\footnotesize, outer sep=1pt}
}

\DeclareMathOperator{\ecd}{ecd}
\DeclareMathOperator{\expand}{expand}

\theoremstyle{definition}
\newtheorem{theorem}{Theorem}
\newtheorem{lemma}[theorem]{Lemma}
\newtheorem{proposition}[theorem]{Proposition}
\newtheorem{corollary}[theorem]{Corollary}

\newtheorem{conjecture}{Conjecture}
\newtheorem{question}[conjecture]{Question}

\title{Even cycle decompositions of small regular graphs}
\author{Alex Diko\\[4pt]
J\'an Maz\'ak\\[4pt]
{\small Comenius University}\\
{\small Bratislava}\\
{\small Slovakia}\\[4pt]
{\small\texttt{jan.mazak@fmph.uniba.sk}}}
\date{30 June 2026}

\begin{document}

\maketitle

\begin{abstract}
An even cycle decomposition (ecd) partitions the edges of a graph into even
cycles. Its ecd index is the minimum $k$ for which the graph has an
ecd whose cycles can be partitioned into $k$ families of pairwise
vertex-disjoint cycles.

We construct simple $4t$-regular graphs with ecd index $3t$ for every integer
$t\geq1$, showing that the gap above half the degree can be arbitrarily
large. For every $k\geq2$, we construct simple $2k$-regular, $2k$-vertex-connected graphs of even order with
ecd index $k+1$. We determine the exact ecd index of the line graph of every connected
simple class 1 regular graph and of large families of circulant graphs.
We also establish theoretical
degree-sequence, cycle-length, and treewidth upper bounds that speed up
our SAT-based algorithms and make exhaustive computations on small
regular graphs feasible.
We explain in a unified fashion the obstructions in previously constructed
4-regular graphs without an ecd and in every such graph in our census.
We propose several conjectures and open questions naturally arising from
our investigation.
\end{abstract}


\paragraph*{MSC-class:} 05C70, 05C38, 05C15


\section{Introduction}

A \emph{cycle} in a graph is a connected $2$-regular subgraph.
A \textit{cycle decomposition} of a graph
is a partition of its edges into cycles such that
each edge belongs to exactly one cycle.
An \textit{even cycle decomposition}, henceforth referred to as \textit{ecd},
is a partition in which each cycle has even length.
All these terms make perfect sense for multigraphs, so unless we specify we are talking about simple graphs,
we allow parallel edges in graphs (but no loops, since a loop cannot be part of an even cycle),
and we allow a cycle in an ecd to consist of two parallel edges.
We will confine the discussion to connected graphs; disconnected graphs can be treated component by component.

A finite connected graph is Eulerian if it has an Euler tour, equivalently
if every vertex has even degree. Splitting an Euler tour at repeated
vertices gives a cycle decomposition. Conversely, each cycle contributes
degree two at each of its vertices, so a cycle decomposition forces all
degrees to be even.
Thus a simple necessary condition for a graph to have an ecd
is that it must be Eulerian. For instance, every Eulerian
bipartite graph has an ecd, since all the cycles in the cycle decomposition of a bipartite graph have even length.

A \emph{block} in a graph is any maximal connected subgraph without a cut vertex.
(If $G$ is $2$-vertex-connected, then the whole $G$ is a single block.)
A cycle of an ecd is confined to a single block. The number of edges in a block is thus
the sum of the lengths of the cycles of an ecd containing edges of that block, leading to the following observation.

\begin{proposition}
If a graph $G$ has an ecd, every block of $G$ has an even number of edges.
\label{prop:block}
\end{proposition}

The question of the existence of ecds has been settled for planar graphs decades ago: every Eulerian planar graph
in which every block has an even number of edges has an ecd \cite{seymour1981even}.
Some Eulerian graphs without cut vertices do not have an ecd:
there are infinitely many such 4-vertex-connected 4-regular graphs \cite{liu20224}.

Given an ecd of a graph $G$, we can colour the cycles
of the ecd so that cycles sharing a vertex belong to different \textit{colour classes}.
If we use $k$ colour classes $1,\ldots,k$ for such a colouring, we
say that $G$ admits a \textit{$k$-ecd}. Some colour classes may be unused
(if a graph $G$ admits a $k$-ecd, it also admits a $(k+1)$-ecd), so we are naturally interested in minimizing $k$.
The \textit{ecd index} of a graph with an ecd is the smallest $k$
such that the graph admits a $k$-ecd.
Ecd index is sometimes also called \textit{ecd size}.
It turns out that it is useful for algorithmic purposes to distinguish edges within each cycle of an ecd: for each non-empty colour class $i$, we can colour the
edges on even positions in the cycles with colour $2i$
and on odd positions with colour $2i-1$.
We call such a colouring a \textit{$k$-ecd colouring}.
Ecd index was first introduced by Bonisoli et al. \cite{bonisoli2017even}.
The inspiration behind our definition of $k$-ecd colouring
is taken from Malnegro et al. \cite{malnegro2024h},
who used a similar characterization for $4$-regular graphs with $2$- and $3$-ecd colourings.
Their paper is focused on $H$-colourings of $4$-regular graphs, a topic which we do not discuss here,
apart from mentioning that having a $3$-ecd is equivalent to being $2C_3$-colourable,
which was our reason to investigate $3$-ecd in more detail.

The following easily proved proposition formally characterizes ecd colourings by local conditions.
Given $f: E(G) \rightarrow \{1, 2, \ldots, 2k\}$, we say that edges $e_1$ and $e_2$ belong
to the same colour class $i$ if $\{f(e_1), f(e_2)\} \subseteq \{2i-1, 2i\}$
for some $i \in \{1, \ldots, k\}$. This characterization is the basis of the SAT encodings in
Section~\ref{sec:algorithms}.

\begin{proposition}
    Let $f: E(G) \rightarrow \{1, 2, \ldots, 2k\}$. Then $f$
    is a $k$-ecd colouring of $G$
    if and only if both of the following conditions hold:
    \begin{itemize}
        \item[(i)] If a vertex $v \in V(G)$ is
        incident to distinct edges $e_1, e_2$
        that belong to the same colour class $i$,
        then $\{f(e_1), f(e_2)\} = \{2i-1, 2i\}$.
        \item[(ii)] If a vertex $v \in V(G)$ is incident to an
        edge $e$ in a colour class $i$, then
        it must be incident to exactly one other
        edge in the same colour class $i$.
    \end{itemize}
    \label{prop:ecd-colouring-characterization}
\end{proposition}

Ecd index can be bounded by common graph invariants such as degree, size, and chromatic index.
Let us write $\ecd(G)$ for the ecd index of $G$; if $G$ does not have an even cycle decomposition, we put $\ecd(G)=\infty$.
The upper bound in the following obvious proposition is improved in Section~\ref{sec:ecd-index-bounds}.

\begin{proposition}
Assume $G$ has an ecd. Then $\ecd(G)\ge \Delta(G)/2$ and $\ecd(G)\le |E|/2$.
\label{prop:elementary-ecd-bounds}
\end{proposition}
For simple graphs with an ecd, the same counting argument gives
$\ecd(G)\leq |E(G)|/4$: every even cycle has at least four edges, and
assigning a distinct colour to each cycle gives a valid ecd colouring.

\begin{proposition}[\cite{bonisoli2017even}]
    A $2k$-regular graph $G$ has ecd index $k$ if and only
    if $\chi'(G) = 2k$.
    Moreover, if $k\geq1$ and $G$ has odd order, then
    $\ecd(G)\geq k+1$.
    \label{prop:ecd-index-chromatic-index}
\end{proposition}
The odd-order assertion follows because a $k$-ecd would make each colour
class a spanning union of even cycles, which is impossible on an odd number
of vertices.
We call a graph $G$ \emph{class 1} if its chromatic index is $\Delta(G)$,
and \emph{class 2} otherwise. The preceding proposition says that an
even-degree regular graph is class 1 if and only if its ecd index is
$\Delta(G)/2$.

\begin{lemma}
\label{lem:bipartite-ecd-index}
Let $B$ be a bipartite multigraph in which every vertex has even degree. Then
\[
    \ecd(B)=\Delta(B)/2.
\]
\end{lemma}
\begin{proof}
Write $\Delta(B)=2k$ and add $(2k-d_B(v))/2$ loops at each vertex to
obtain a $2k$-regular multigraph $B'$. Orient the Eulerian components of
$B'$ so that every vertex has indegree and outdegree $k$, and split each
vertex into an out-copy and an in-copy. Each directed edge joins the
out-copy of its tail to the in-copy of its head, so the edges between
the two copies of the same vertex correspond exactly to the added loops.
The resulting $k$-regular bipartite multigraph decomposes into $k$ perfect
matchings. Identifying the two copies of each vertex again, each matching
gives a spanning $2$-regular subgraph of $B'$, since it meets each copy
once. Each component is an added loop or an even cycle of $B$, since
$B$ is bipartite. Deleting the loops therefore gives a $k$-ecd of $B$,
and the degree lower bound of Proposition~\ref{prop:elementary-ecd-bounds}
gives equality. Casselgren et al.~\cite{casselgren2025near} used a similar
loop-padding argument with Petersen's 2-factor theorem.
\end{proof}

One of the first explorations of ecds was conducted by Markstr\"om
\cite{MARKSTROM20122676}. He mentioned an infinite class of $4$-regular graphs without an ecd: cut an edge in $K_5$ and insert this gadget into any edge of a 4-regular graph.
His results were later improved by the construction of infinite families
of 3-vertex-connected \cite{mavcajova2013even} and 4-vertex-connected
\cite{liu20224} simple $4$-regular graphs without an ecd, which we
revisit in Section~\ref{sec:three-connected-no-ecd}.
For every fixed $d\geq3$, almost all $d$-regular graphs of even order are
class 1 \cite{kim2001random}. Hence, by
Proposition~\ref{prop:ecd-index-chromatic-index}, almost all $2k$-regular
graphs of even order have ecd index $k$. For 4-regular graphs of odd
order, Markstr\"om conjectured that a random graph asymptotically almost
surely decomposes into a cycle of length one less than its order and two
even cycles, and therefore has a 3-ecd
\cite[Conjecture~2.3]{MARKSTROM20122676}.

The problem of determining whether a graph has an ecd appears to be hard both algorithmically and theoretically;
no simple general characterization of graphs with an ecd is known. We will present various open questions throughout the paper,
interspersed with related computational observations and proofs of ecd-related properties for specific families of graphs.
The present paper develops and extends results from Diko's bachelor thesis
\cite{diko2026thesis}. Its census of small regular graphs in
Section~\ref{sec:computational-results} motivates a connectivity-based
analysis of the 4-regular case in Sections~\ref{sec:connectivity-1}--
\ref{sec:three-connected-no-ecd}.

We use the following terminology for multipoles, following
\cite[Chapter~1]{diko2026thesis}. A \emph{multipole} is a
multigraph that may have \emph{dangling edges}, each of which has exactly
one end incident with a vertex; the other end is free. A \emph{$k$-pole}
is a multipole with $k$ dangling edges. Splitting
an edge replaces it by two dangling edges, one incident with each of its
endvertices, while deleting a vertex turns each edge incident with it into
a dangling edge. Joining two dangling edges removes them and adds an edge
between their incident vertices.


\section{Upper bounds on ecd index}
\label{sec:ecd-index-bounds}

The upper bound in Proposition~\ref{prop:elementary-ecd-bounds} can be
improved substantially. Let $\mathcal D$ be an ecd of $G$.
The \emph{conflict graph} $H_{\mathcal D}$
has one vertex for every cycle of $\mathcal D$, with two vertices adjacent
when the corresponding cycles meet in $G$. The minimum number of colours
needed by $\mathcal D$ is $\chi(H_{\mathcal D})$, and hence
\[
    \ecd(G)=\min_{\mathcal D}\chi(H_{\mathcal D}),
\]
where the minimum is over all ecds of $G$.
For a graph $G$ with even degrees, put
\[
    S(G)=\sum_{v\in V(G)}
        \frac{d_G(v)}{2}\left(\frac{d_G(v)}{2}-1\right),
\]
and let $g_e(G)$ denote the length of a shortest even cycle of $G$.

\begin{lemma}[Critical conflict counting]
\label{lem:conflict-counting}
Let $\mathcal D$ be a cycle decomposition of a graph $G$ (so all
degrees are even), and set $q=\chi(H_{\mathcal D})$. Then
$H_{\mathcal D}$ has a subgraph $J$ with $|V(J)|\geq q$ such that every
cycle $C$ represented in $J$ satisfies
\begin{equation}
    q-1\leq d_J(C)
       \leq \sum_{v\in V(C)}\left(\frac{d_G(v)}{2}-1\right).
    \label{eq:critical-conflict-core}
\end{equation}
Consequently,
\begin{equation}
    q(q-1)\leq S(G).
    \label{eq:degree-sequence-ecd-bound}
\end{equation}
\end{lemma}

\begin{proof}
Starting from $H_{\mathcal D}$, repeatedly delete a vertex whose removal
leaves the current subgraph with chromatic number $q$. The
resulting subgraph $J$ is vertex-critical and $q$-chromatic. Hence
$|V(J)|\geq q$ and $J$ has minimum degree at least $q-1$: otherwise a $(q-1)$-colouring of
$J-x$ could be extended to a vertex $x$ of degree at most $q-2$.

Let $C$ be a cycle represented in $J$. Every neighbour of $C$ in $J$ meets
$C$ at some vertex $v$. Besides $C$, exactly $d_G(v)/2-1$ cycles of
$\mathcal D$ contain $v$. Summing over the vertices of $C$ may count a
neighbour more than once. Therefore~\eqref{eq:critical-conflict-core} holds.
Summing~\eqref{eq:critical-conflict-core} over the cycles represented in $J$ gives
\[
    q(q-1)\leq\sum_{C\in V(J)}d_J(C)
    \leq S(G),
\]
since each vertex $v$ belongs to at most $d_G(v)/2$ represented cycles.
\end{proof}

\begin{theorem}
\label{thm:conflict-ecd-bound}
Let $G$ be a multigraph with an ecd and $\Delta(G)\geq4$. Then
\begin{equation}
    \ecd(G)\leq
    \min\left\{
        \left\lfloor
            \frac{1+\sqrt{1+4S(G)}}{2}
        \right\rfloor,
        \left\lfloor |E|/g_e(G)\right\rfloor,
        \max\left\{j\geq1:
            2j\left\lceil\frac{j-1}{\Delta(G)-2}\right\rceil\leq |E|
        \right\}
    \right\}.
    \label{eq:square-root-ecd-bound}
\end{equation}
In particular, if $G$ is 4-regular of order $n$ and $\ecd(G)=k$, then
$k\lceil(k-1)/2\rceil\leq n$.
\end{theorem}

\begin{proof}
Let $k=\ecd(G)$ and choose an ecd $\mathcal D$ with
$\chi(H_{\mathcal D})=k$. Apply Lemma~\ref{lem:conflict-counting} to obtain
$J$ and the inequality $k^2-k-S(G)\leq0$. Solving this quadratic inequality gives
the first bound in~\eqref{eq:square-root-ecd-bound}.

The subgraph $J$ represents at least $k$ edge-disjoint even cycles, each
of length at least $g_e(G)$. Thus $|E|\geq k g_e(G)$, giving the second bound.
For every cycle $C$ represented in $J$, inequality
\eqref{eq:critical-conflict-core} gives
\[
    k-1\leq\sum_{v\in V(C)}\left(\frac{d_G(v)}{2}-1\right)
       \leq \frac{\Delta(G)-2}{2}|E(C)|.
\]
Since $|E(C)|$ is even, each of these cycles has length at least
$2\lceil(k-1)/(\Delta(G)-2)\rceil$. Summing over the cycles gives
\begin{equation}
    |E|\geq 2k\left\lceil\frac{k-1}{\Delta(G)-2}\right\rceil.
    \label{eq:parity-ecd-bound}
\end{equation}
Hence $j=k$ belongs to the set in the third bound.
For a 4-regular graph of order $n$, substituting $|E|=2n$ and $\Delta(G)=4$
into~\eqref{eq:parity-ecd-bound} gives the asserted inequality.
\end{proof}

For example, if $G$ is 4-regular of order 19 or 20 and has an
ecd, then Theorem~\ref{thm:conflict-ecd-bound} gives $\ecd(G)\leq6$.

The bound in Theorem~\ref{thm:conflict-ecd-bound} is tight for
multigraphs: on three vertices with $2t$ parallel edges joining each pair,
where $t\geq1$, every ecd consists of $3t$ pairwise intersecting 2-cycles,
so the index is $3t=|E|/g_e$.
For simple graphs, the construction in
Theorem~\ref{thm:simple-regular-large-index} gives $4t$-regular examples
of order $48t^2+3$ and index $3t=3\Delta/4$.
For these graphs, the upper bound in
Theorem~\ref{thm:conflict-ecd-bound} grows quadratically with $t$, so its
ratio to the actual index grows linearly.
For 4-regular graphs, the upper bound in
Theorem~\ref{thm:conflict-ecd-bound} is asymptotic to $\sqrt{2n}$,
whereas our examples attain index 3. The gap between the bound of
Theorem~\ref{thm:conflict-ecd-bound} and the actual ecd index can thus be
large for simple graphs, yet simple regular graphs with larger finite index
are hard to find in small degrees (see Section~\ref{sec:computational-results}).
Tightness for the multigraph family above suggests that improving the bound
in general may not be easy. We consider this gap an important open question.

A tree decomposition gives a bound on ecd index independent of the order. We use the
standard notation $\operatorname{tw}(G)$ for treewidth.

\begin{theorem}[Tree-decomposition bound]
\label{thm:tree-decomposition-ecd-bound}
Let $G$ have an ecd. Then, taking the minimum over all tree decompositions
$(T,\{B_x:x\in V(T)\})$ of $G$, we have
\[
    \ecd(G)\leq
    \min_{(T,\{B_x\})}
    \max_{x\in V(T)}\sum_{v\in B_x}\frac{d_G(v)}{2}
    \leq
    \frac{(\operatorname{tw}(G)+1)\Delta(G)}{2}.
\]
\end{theorem}

\begin{proof}
Fix a tree decomposition $(T,\{B_x:x\in V(T)\})$ of $G$ and set
\[
    U=\max_{x\in V(T)}\sum_{v\in B_x}\frac{d_G(v)}{2}.
\]
Fix an ecd $\mathcal D$ of $G$. For each $x\in V(T)$, let $B'_x$ be the
set of cycles of $\mathcal D$ that contain at least one vertex of $B_x$.
These sets form the bags of a tree decomposition of the conflict graph
$H_{\mathcal D}$. Each cycle appears in a derived bag because any vertex
of the cycle belongs to an original bag. Intersecting cycles are covered by any bag
containing a common vertex. Moreover, if $T_v$ is the subtree formed by
the bags containing $v$, then the bags whose derived sets contain a fixed
cycle $C$ are indexed by
\[
    \bigcup_{v\in V(C)}V(T_v).
\]
This union is connected because $T_u\cap T_v\neq\emptyset$ whenever
$uv\in E(C)$.
Exactly $d_G(v)/2$ cycles of $\mathcal D$ contain $v$. Therefore
\[
    |B'_x|\leq\sum_{v\in B_x}\frac{d_G(v)}{2}\leq U
\]
for every $x\in V(T)$.
If $U=0$, every derived bag is empty, so $\mathcal D$ is empty and
$\ecd(G)=0$. Otherwise $\operatorname{tw}(H_{\mathcal D})\leq U-1$ and the standard inequality
$\chi(H)\leq\operatorname{tw}(H)+1$~\cite[Section~1.1, p.~3]{brianski2026defective}
gives $\ecd(G)\leq\chi(H_{\mathcal D})\leq U$.

Taking the minimum over all decompositions proves the first inequality. In
a decomposition of width $\operatorname{tw}(G)$, every bag has at most
$\operatorname{tw}(G)+1$ vertices and $d_G(v)/2\leq\Delta(G)/2$, which proves
the second.
\end{proof}

\phantomsection\label{disc:planar-ecd-bounds}
Every outerplanar graph has treewidth at most
2~\cite[Section~3]{bodlaender1993tourist}, so
Theorem~\ref{thm:tree-decomposition-ecd-bound} gives
$\ecd(G)\leq 3\Delta(G)/2$ for every outerplanar graph $G$ admitting an ecd.
Seymour's theorem~\cite[Theorem~3.1]{seymour1981even} gives existence of an
ecd for Eulerian planar graphs with an even number of edges in each block.
His proof does not appear to give better bounds on the ecd index of planar
or outerplanar graphs than Theorems~\ref{thm:conflict-ecd-bound}
and~\ref{thm:tree-decomposition-ecd-bound}. It is natural to ask whether any
bound of the form $C\Delta(G)$, for an absolute constant $C$, holds for all
planar graphs admitting an ecd. For degree 4, Malnegro and Ozeki
ask whether every 2-connected planar 4-regular graph has ecd index at most
3~\cite[Section~4.3]{malnegro2024h}.


\section{Algorithms}
\label{sec:algorithms}

For every fixed $d\geq3$, deciding whether a simple $d$-regular graph is
$d$-edge-colourable is NP-complete \cite{leven1983np}. Hence, by
Proposition~\ref{prop:ecd-index-chromatic-index}, for every fixed $k\geq2$,
deciding whether a simple $2k$-regular graph admits a $k$-ecd is
NP-complete, deciding whether it satisfies $\ecd(G)>k$ is co-NP-complete, and
computing its ecd index is NP-hard.

These complexity results do not settle the unrestricted
problem of deciding whether a graph has some ecd. This problem belongs to NP,
since an ecd is a polynomial-size certificate. As far as we know, it is not
known whether this unrestricted problem is NP-hard.

In this section, we present and compare the most efficient algorithms we have
developed for finding even cycle decompositions; a more detailed discussion
is given in the thesis \cite{diko2026thesis}.

\subsection{SAT and backtracking}

Let $G=(V,E)$. For $X\subseteq V$, let $\delta(X)$ denote the set of edges
with exactly one end in $X$, and write $\delta(v)=\delta(\{v\})$ for the set
of edges incident with $v$. Both encodings below use, for every $e\in E$, a
Boolean variable $p_e$, which is true if and only if $e$ receives the even
colour of its class.

\bigskip

\noindent\textbf{Ecd existence without a size bound.}
An ecd partitions the edges into colour classes. For every unordered pair of
distinct edges $e,f\in E$, introduce a Boolean variable $s_{e,f}$, which is
true if and only if $e$ and $f$ belong to the same colour class.

By Proposition \ref{prop:ecd-colouring-characterization}, $G$ has an ecd if and only
if the following set of clauses is satisfiable:
\begin{enumerate}
    \item Two incident edges in the same colour class have opposite parities.
    \[
        \lnot s_{e,f}\lor\lnot p_e\lor\lnot p_f
        \qquad\text{and}\qquad
        \lnot s_{e,f}\lor p_e\lor p_f
        \qquad
        \begin{gathered}
            \text{for every }v\in V\text{ and every unordered}\\
            \text{pair of distinct edges }e,f\in \delta(v).
        \end{gathered}
    \]

    \item Each endpoint of an edge is incident with another edge in the same
    colour class.
    \[
        \bigvee_{f\in \delta(v)\setminus\{e\}}s_{e,f}
        \qquad\text{for every }v\in V\text{ and }e\in \delta(v).
    \]

    \item Belonging to the same colour class is transitive. For every
    unordered triple of distinct edges $e,f,g\in E$, include
    \[
        \lnot s_{e,f}\lor\lnot s_{f,g}\lor s_{e,g}
    \]
    and the two clauses obtained by cyclically permuting $e,f,g$.
\end{enumerate}

The encoding uses $O(|E|^2)$ variables and
$O(|E|\Delta(G)+|E|^3)$ clauses.
To break the symmetry obtained by swapping the two colours in a class, fix a
total order $<$ on $E$ and add the following clauses:
\begin{enumerate}
    \setcounter{enumi}{3}
    \item The smallest edge in every colour class has even parity.
    \[
        p_e\lor\bigvee_{\substack{f\in E\\f<e}}s_{e,f}
        \qquad\text{for every }e\in E.
    \]
\end{enumerate}
These constraints add $O(|E|)$ clauses and no variables, without changing
satisfiability.

The transitivity clauses cannot be omitted or restricted to adjacent edges:
local-only clauses admit assignments which do not represent ecds. A
counterexample and a detailed explanation are given in Diko's thesis
\cite[Chapter~2, Section~2.1]{diko2026thesis}.

\bigskip

\noindent\textbf{Finding a $k$-ecd.}
The following clauses and variables encode a Boolean
formula in CNF that is satisfiable if and only if $G$
admits a $k$-ecd (by Proposition \ref{prop:ecd-colouring-characterization}).
For every edge $e\in E$ and every colour class $i\in\{1,\ldots,k\}$,
introduce a Boolean variable $c_{e,i}$, which is true if and only if
$e$ belongs to class $i$; the variables $p_e$ are as above.
\begin{enumerate}
    \item Every edge belongs to at least one colour class.
    \[
        \bigvee_{i=1}^{k}c_{e,i}
        \qquad\text{for every }e\in E.
    \]

    \item Every edge belongs to at most one colour class.
    \[
        \lnot c_{e,i}\lor\lnot c_{e,j}
        \qquad\text{for every }e\in E\text{ and }1\leq i<j\leq k.
    \]

    \item Two incident edges in the same colour class have opposite parities.
    \[
        \begin{aligned}
        &\lnot c_{e,i}\lor\lnot c_{f,i}\lor\lnot p_e\lor\lnot p_f,\\
        &\lnot c_{e,i}\lor\lnot c_{f,i}\lor p_e\lor p_f
        \end{aligned}
    \]
    for every $v\in V$, every unordered pair of distinct edges
    $e,f\in \delta(v)$, and every $i\in\{1,\ldots,k\}$.

    \item Each endpoint of an edge is incident with another edge in the same
    colour class.
    \[
        \lnot c_{e,i}\lor
        \bigvee_{f\in \delta(v)\setminus\{e\}}c_{f,i}
        \qquad
        \begin{gathered}
            \text{for every }v\in V,\ e\in \delta(v),\\
            \text{and }i\in\{1,\ldots,k\}.
        \end{gathered}
    \]
\end{enumerate}

Although clause 4 explicitly requires only another edge of class $i$ at
an endpoint, clause 3 makes this edge unique: among three incident edges
of one class, two would have the same parity.
The encoding uses $(k+1)|E|$ variables and
$O(k^2|E|+k|E|\Delta(G))$ clauses.

Permuting the colour classes or swapping the two colours within a class
produces equivalent satisfying assignments. To break these symmetries,
fix a total order $<$ on $E$ and add the following clauses:

\begin{enumerate}
    \setcounter{enumi}{4}
    \item The smallest edge in each colour class has even parity.
    \[
        \lnot c_{e,i}\lor p_e\lor
        \bigvee_{\substack{f\in E\\f<e}}c_{f,i}
        \qquad
        \begin{gathered}
            \text{for every }e\in E,\\
            \text{and }i\in\{1,\ldots,k\}.
        \end{gathered}
    \]

    \item Before an edge belongs to a colour class other than the first, a
    smaller edge belongs to the preceding class.
    \[
        \lnot c_{e,i}\lor
        \bigvee_{\substack{f\in E\\f<e}}c_{f,i-1}
        \qquad
        \begin{gathered}
            \text{for every }e\in E,\\
            \text{and }i\in\{2,\ldots,k\}.
        \end{gathered}
    \]
\end{enumerate}

Here an empty disjunction is false; in particular, these clauses force the
smallest edge to belong to class $1$ and have even parity.
These constraints add $O(k|E|)$ clauses and no variables, without changing
satisfiability.

\paragraph{Backtracking.}
We also implemented an independent backtracking algorithm based on the same
local characterization. It assigns colours directly to edges, incorporates
symmetry breaking, and prunes partial assignments that cannot be completed.
Besides providing an independent check of the SAT implementation, it can
determine ecd index and produce a colouring in a single run. In contrast,
the bounded SAT encoding may require several calls with different values of
$k$, so the backtracking algorithm is not a priori inferior.

\subsection{Comparison of running times}

We benchmarked simple 4-regular class 2 graphs of order 17 and vertex
connectivity at least 3. The test set comprised a fixed pseudorandom sample
of 100 graphs of ecd index 3 and all 137 graphs in our census without an ecd.
We checked these classifications using backtracking and both SAT encodings,
with bounds 3 and 5 for the bounded encoding. The timings were taken on an
AMD Ryzen 7 8845HS and include formula construction, symmetry breaking,
solving, and validation of satisfying assignments. We used
Clasp~\cite{clasp}, which was faster than Kissat in Diko's solver
comparison~\cite[Chapter~2, Section~4]{diko2026thesis}.

In the column ``Index 3'', Table~\ref{tab:ecd-benchmark-index} compares
exact-index backtracking with a single SAT test for a $3$-ecd. Class 2
already rules out a $2$-ecd by Proposition~\ref{prop:ecd-index-chromatic-index}.
In the census for this class, only 137 graphs fail the $3$-ecd test, and
none has an ecd (Table~\ref{tab:ecd-4-regular}). The last two columns
compare a $5$-ecd test with unrestricted SAT on these graphs.
Theorem~\ref{thm:conflict-ecd-bound} makes the bound $5$ sufficient.
The first row gives the unrestricted backtracking times.

\begin{table}[H]
\centering
\small
\begin{tabular}{lr@{ }lr@{ }lr@{ }l}
\hline
\textbf{Method} & \multicolumn{2}{c}{\textbf{Index 3}}
    & \multicolumn{2}{c}{\textbf{No ecd, $k=5$}}
    & \multicolumn{2}{c}{\textbf{No ecd, unrestricted}} \\
\hline
Backtracking, exact index & 118.2 & [106.7--132.3]
    & \multicolumn{2}{c}{---} & 64 & [60--73] \\
\hline
SAT, no sym. breaking & 0.4 & [0.4--0.4]
    & 847 & [658--1101] & 39 & [35--46] \\
SAT, explicit & 0.5 & [0.5--0.6]
    & 40 & [34--46] & 26 & [23--29] \\
SAT, BreakID & 2.5 & [2.4--2.6]
    & 62 & [52--69] & 57 & [52--64] \\
\hline
\end{tabular}
\caption{Running times in milliseconds, as ``median [interquartile range]''.}
\label{tab:ecd-benchmark-index}
\end{table}

On the index-3 graphs, all three SAT variants had lower median times than
backtracking, with the lowest median for the test without symmetry breaking. On the
negative set, explicit symmetry breaking gave the lowest median for both
SAT encodings, and unrestricted SAT was fastest despite its larger formula.
BreakID improved only the $5$-ecd timings on the negative set.

To compare the methods on larger graphs, we tested for an ecd in the line
graphs of a small sample of 30-vertex snarks, as in the computational
checks of Markstr\"om's line-graph conjecture (Section~\ref{sec:line-graphs}).
These line graphs have 45 vertices and 90 edges.
For the unbounded existence question, Theorem~\ref{thm:conflict-ecd-bound} shows that it is enough to use
bound 9 (the third bound from the theorem).
The elementary bound $|E|/4$ for simple graphs gives bound 22
for these line graphs.
The SAT rows in Table~\ref{tab:ecd-benchmark-line-graphs} report individual
Clasp tests with explicit symmetry breaking on these line graphs and on the
order-17 index-3 sample. Diko's thesis reports similar timing comparisons
on other test sets~\cite[Chapter~2, Section~4]{diko2026thesis}.

\begin{table}[H]
\centering
\small
\begin{tabular}{lr@{ }lr@{ }l}
\hline
\textbf{Method} & \multicolumn{2}{c}{\textbf{Order 17, index 3}}
    & \multicolumn{2}{c}{\textbf{Line graphs of order-30 snarks}} \\
\hline
SAT, bound 3 & 0.5 & [0.5--0.6] & 2.4 & [2.3--2.4] \\
SAT, Thm.~\ref{thm:conflict-ecd-bound} bound
    & 0.9 & [0.9--1.0] ($5$-ecd) & 11.5 & [9.4--12.6] ($9$-ecd) \\
SAT, $|E|/4$ bound
    & 1.7 & [1.6--1.8] ($8$-ecd) & 29.1 & [27.3--33.8] ($22$-ecd) \\
SAT, unrestricted & 7.2 & [7.1--7.6] & 148.7 & [141.7--168.8] \\
Backtracking, exact index
    & 118.2 & [106.7--132.3] & $>60\,000$ & \\
\hline
\end{tabular}
\caption{Scaling benchmark in milliseconds, as ``median [interquartile range]''.}
\label{tab:ecd-benchmark-line-graphs}
\end{table}

\section{Computational results}
\label{sec:computational-results}

We generated nonisomorphic connected regular graphs using
\texttt{genreg}~\cite{meringer1999fast} in graph6 format for simple graphs
and \texttt{nauty}~\cite{mckay2014practical} in sparse6 format for
multigraphs. Tables~\ref{tab:ecd-4-regular}
and~\ref{tab:ecd-4-regular-multigraphs} give the ecd census for class-2
4-regular graphs by vertex connectivity $\kappa$ and edge connectivity
$\lambda$. The simple census goes through order 20 and the non-simple
multigraph census through order 15.
The thesis splits its counts of simple graphs without an ecd at
$(\kappa,\lambda)=(2,2)$ by whether they contain the $2$-pole obtained by
splitting an edge of $K_5$~\cite[Chapter~3, Section~3.3]{diko2026thesis}.
Every class-2 4-regular graph in our census either admits a $3$-ecd or has
no ecd.

\begin{table}[!t]
\centering
\footnotesize
\setlength{\tabcolsep}{4pt}
\textbf{(a) Graphs with ecd index 3}\par\smallskip
\begin{tabular}{r|rrrr|r}
\hline
& \multicolumn{4}{c|}{Connectivity $(\kappa,\lambda)$} & \\
\textbf{Order} & $(2,2)$ & $(2,4)$ & $(3,4)$ & $(4,4)$ & \textbf{Total} \\
\hline
5  & 0 & 0 & 0 & 0 & 0 \\
6  & 0 & 0 & 0 & 0 & 0 \\
7  & 0 & 0 & 1 & 1 & 2 \\
8  & 0 & 0 & 0 & 0 & 0 \\
9  & 0 & 0 & 3 & 12 & 15 \\
10 & 0 & 0 & 0 & 0 & 0 \\
11 & 0 & 3 & 55 & 204 & 262 \\
12 & 0 & 0 & 0 & 0 & 0 \\
13 & 14 & 67 & 1562 & 9088 & 10731 \\
14 & 29 & 0 & 0 & 0 & 29 \\
15 & 445 & 2075 & 80188 & 721845 & 804553 \\
16 & 910 & 0 & 0 & 0 & 910 \\
17 & 18891 & 100430 & 6055208 & 80003733 & 86178262 \\
18 & 46651 & 0 & 1 & 0 & 46652 \\
19 & \multicolumn{4}{c|}{not separated by connectivity} & 11943280611 \\
20 & 3233306 & 1 & 13 & 0 & 3233320 \\
\hline
\end{tabular}

\medskip
\textbf{(b) Graphs without an ecd}\par\smallskip
\begin{tabular}{r|rrrr|r}
\hline
& \multicolumn{4}{c|}{Connectivity $(\kappa,\lambda)$} & \\
\textbf{Order} & $(1,2)$ & $(2,2)$ & $(3,4)$ & $(4,4)$ & \textbf{Total} \\
\hline
5  & 0 & 0 & 0 & 1 & 1 \\
6  & 0 & 0 & 0 & 0 & 0 \\
7  & 0 & 0 & 0 & 0 & 0 \\
8  & 0 & 0 & 0 & 0 & 0 \\
9  & 0 & 0 & 1 & 0 & 1 \\
10 & 0 & 1 & 0 & 0 & 1 \\
11 & 1 & 2 & 0 & 0 & 3 \\
12 & 2 & 7 & 0 & 0 & 9 \\
13 & 10 & 28 & 9 & 0 & 47 \\
14 & 42 & 126 & 0 & 0 & 168 \\
15 & 210 & 728 & 0 & 0 & 938 \\
16 & 1167 & 4749 & 0 & 0 & 5916 \\
17 & 7445 & 35790 & 136 & 1 & 43372 \\
18 & 53662 & 300855 & 0 & 0 & 354517 \\
19 & 432770 & 2774266 & 0 & 0 & 3207036 \\
20 & 3843536 & 27639839 & 0 & 0 & 31483375 \\
\hline
\end{tabular}
\caption[Ecd properties of simple 4-regular graphs]{All simple 4-regular
graphs of orders 5--20 that are not 4-edge-colourable (equivalently, do not
admit a 2-ecd), classified by ecd index and connectivity.}
\label{tab:ecd-4-regular}
\end{table}

\begin{table}[!t]
\centering
\scriptsize
\setlength{\tabcolsep}{4pt}
\begin{tabular}{r|rrr|r|rr|r}
\hline
& \multicolumn{4}{c|}{Ecd index 3} & \multicolumn{3}{c}{No ecd} \\
\textbf{Order}
& $(2,2)$ & $(2,4)$ & $(3,4)$ & \textbf{Total}
& $(1,2)$ & $(2,2)$ & \textbf{Total} \\
\hline
3  & 0 & 1 & 0 & 1 & 0 & 0 & 0 \\
4  & 0 & 0 & 0 & 0 & 0 & 0 & 0 \\
5  & 1 & 2 & 1 & 4 & 1 & 0 & 1 \\
6  & 1 & 0 & 0 & 1 & 2 & 0 & 2 \\
7  & 14 & 14 & 11 & 39 & 8 & 1 & 9 \\
8  & 11 & 0 & 0 & 11 & 27 & 1 & 28 \\
9  & 266 & 204 & 223 & 693 & 116 & 7 & 123 \\
10 & 279 & 0 & 1 & 280 & 499 & 13 & 512 \\
11 & 7987 & 5999 & 8399 & 22385 & 2502 & 75 & 2577 \\
12 & 9447 & 3 & 2 & 9452 & 13722 & 271 & 13993 \\
13 & 377959 & 305013 & 523374 & 1206346 & 85080 & 1521 & 86601 \\
14 & 463441 & 73 & 58 & 463572 & 584301 & 8299 & 592600 \\
15 & 26063985 & 22810259 & 46050997 & 94925241 & 4444575 & 54810 & 4499385 \\
\hline
\end{tabular}
\caption[Ecd properties of 4-regular multigraphs]{All non-simple
4-regular multigraphs of orders 3--15 that are not 4-edge-colourable
(equivalently, do not admit a 2-ecd), classified by ecd outcome and
connectivity. Column headings give $(\kappa,\lambda)$.}
\label{tab:ecd-4-regular-multigraphs}
\end{table}

The simple-graph census covers fewer orders in degrees 6 and 8.
For degree 6, odd order gives an odd number of edges and precludes an ecd.
Orders 8, 10, and 12 contain only class-1 graphs. The seven class-2 graphs
of order 14 each have an edge cut of size less than 6 and index 4.

For degree 8, all graphs of orders 9, 11, and 13 have index 5, and all
graphs of orders 10, 12, and 14 are class 1~\cite[Chapter~3]{diko2026thesis}.

In the simple-graph census, every 6- or 8-regular graph with an even
number of edges admits an ecd, and every 4- or 6-regular graph admitting
an ecd has index at most one above half its degree.

\begin{question}
    \label{qu:two-connected-no-ecd}
    Does there exist a simple 2-vertex-connected 6- or 8-regular graph
    with an even number of edges and no ecd?
\end{question}

\begin{question}
    \label{qu:simple-index-above-k-plus-one}
    Does there exist a simple 4- or 6-regular graph admitting an ecd
    whose ecd index is at least 4 or 5, respectively?
\end{question}

\clearpage
\section{Graphs with vertex connectivity 1}
\label{sec:connectivity-1}

We use the following elementary parity obstruction, proved e.g. in
\cite[Chapter~3, Section~2]{diko2026thesis}.
\begin{proposition}
    Let $G$ be a $d$-regular graph, where $d\geq1$. If there is a set
    $S\subseteq V(G)$ of odd cardinality with $|\delta(S)|<d$,
    then $G$ is not $d$-edge-colourable.
    \label{prop:edge-cut-no-d-edge-colouring}
\end{proposition}
\begin{corollary}
    \label{cor:cut-vertex-no-k-ecd}
    A $d$-regular graph with a cut vertex is not $d$-edge-colourable.
    In particular, a $2k$-regular graph with a cut vertex does not admit a $k$-ecd.
\end{corollary}
\begin{proof}
Let $v$ be a cut vertex of $G$. If $G$ has odd order, set $S=V(G)$.
Otherwise, let $S$ be the vertex set of an odd-order component of $G-v$.
In the latter case $|\delta(S)|<d$ since $v$ has neighbours in another
component. In both cases Proposition~\ref{prop:edge-cut-no-d-edge-colouring}
rules out a $d$-edge-colouring. The ecd assertion follows from
Proposition~\ref{prop:ecd-index-chromatic-index}.
\end{proof}

For $4$-regular graphs, the stronger conclusion below
is proved in Diko's thesis~\cite[Chapter~3, Section~2]{diko2026thesis}.
Bonvicini and Mazzuoccolo also note that a cut vertex forces an odd endblock
\cite[Section~2.2]{bonvicini2016edge}.
\begin{proposition}[Cut-vertex obstruction]
    Let $G$ be a 4-regular graph with a cut vertex. Then $G$ does not have an ecd.
    \label{prop:4reg-cut-vertex-no-ecd}
\end{proposition}
For every $k\geq3$ we construct two simple $2k$-regular graphs with a cut
vertex, one with no ecd and one with ecd index $k+1$. In each construction
we subdivide a total of $k$ edges in disjoint copies of $2k$-regular
graphs and identify all the subdivision vertices. The edges chosen in
each copy are vertex-disjoint, so the resulting graph is simple. The
common vertex has degree $2k$ and separates the copies.

For the graph without an ecd, take two copies of $K_{2k+2}$ minus a
perfect matching and $k-2$ copies of $K_{2k+1}$, subdividing one edge in
each. The first two modified copies have $k(2k+2)+1$ edges each, and the
remaining ones have $k(2k+1)+1$ each. Thus exactly two copies have odd
size when $k$ is odd and all $k$ do when $k$ is even. The whole graph has
an even number of edges. Every cycle is confined to one copy, so an
odd-size copy rules out an ecd. In degrees 6 and 8, these graphs show
why Question~\ref{qu:two-connected-no-ecd} requires 2-vertex-connectivity.

For the graph of index $k+1$, take $k$ copies of $K_{2k+1}$ and subdivide
one edge in each if $k$ is odd. If $k$ is even, take $k/2$ copies and
subdivide two vertex-disjoint edges in each. This extends the degree-6
and degree-8 examples in Diko's thesis~\cite[Chapter~3, Section~2]{diko2026thesis}.
By the solution of Alspach's conjecture~\cite{bryant2014cycle}, $K_{2k+1}$
decomposes into $k$ cycles of length $2k$ and one of length $k$ for odd
$k$, or into $k-1$ cycles of length $2k$ and two of lengths $k+1$ and
$2k-1$ for even $k$. In the even case the two odd cycles have length at
least five, so we can pick a vertex-disjoint pair of edges, one on each.
By symmetry, we may therefore arrange that each subdivided edge lies on
a different odd cycle. This gives $k+1$ even cycles per copy. Give them
distinct colours within each copy and permute the colours so that the
$k$ cycles through the cut vertex receive distinct colours. The resulting
$(k+1)$-ecd is optimal by Corollary~\ref{cor:cut-vertex-no-k-ecd}.

\section{Graphs with vertex connectivity 2}
\label{sec:connectivity-2}

Every simple 4-regular graph of vertex connectivity 2 without an ecd in
Table~\ref{tab:ecd-4-regular}(b) has a 2-edge cut.
If exactly two edges join the parts $A$ and $B$ of a partition of $V(G)$,
those edges form a 2-edge cut. The cut is \emph{proper} if its two edges
have distinct ends in each part. In a 2-vertex-connected 4-regular graph
every 2-edge cut is proper, since a common end would be a cut vertex. For
a proper cut, the \emph{closures} $G_A$ and $G_B$ are obtained from $G[A]$
and $G[B]$, respectively, by adding an edge between the two cut ends in
each part. The added edges may be parallel to existing edges.

\begin{theorem}[Decomposition at a 2-edge cut]
\label{thm:two-edge-cut-decomposition}
Let $A$ and $B$ be the parts separated by a proper 2-edge cut in a
multigraph $G$, with closures $G_A$ and $G_B$. If at least one closure has
an even number of edges, then for every $k\geq1$, $G$ admits a $k$-ecd
if and only if both $G_A$ and $G_B$ do.
\end{theorem}

\begin{proof}
By interchanging the parts if necessary, assume that $G_A$ has an even
number of edges, so $G[A]$ has an odd number of edges.
Write $x_1y_1$ and $x_2y_2$ for the cut edges, with $x_i\in A$ and
$y_i\in B$. Let $\mathcal D$ be a $k$-ecd of $G$. Each cycle crosses the cut
an even number of times, so both cut edges lie on the same cycle $C$ of
$\mathcal D$. Deleting them from $C$
leaves an $x_1x_2$-path $P_A$ and a $y_1y_2$-path $P_B$. The path
$P_A$, together with the cycles of $\mathcal D$ contained in $G[A]$,
partitions $E(G[A])$. Since $A$ spans an odd number of edges, $P_A$ is odd.
The cycle $C$ is even, so $P_B$ is odd as well. Replacing $C$ by
$P_A+x_1x_2$ and by $P_B+y_1y_2$ gives ecds of $G_A$ and $G_B$.
Keep the colours of the retained cycles and give each new cycle the
colour of $C$ to obtain a $k$-ecd of each closure.

Conversely, take a $k$-ecd of each closure. Deleting the added edge from
the cycle that contains it leaves an odd path in each part. Joining
these paths by $x_1y_1$ and $x_2y_2$ gives one even cycle of $G$; all
other cycles are unchanged. Permute the colours in one closure so that the
two cycles containing the added edges have the same colour. Give the
joined cycle that colour to obtain a $k$-ecd of $G$.
\end{proof}

In terms of the index, Theorem~\ref{thm:two-edge-cut-decomposition} gives
$\ecd(G)=\max\{\ecd(G_A),\ecd(G_B)\}$. If every vertex degree of $G$ is
divisible by four, then
$\sum_{v\in A}d_G(v)=2|E(G[A])|+2$ implies that $A$ spans an odd number of edges.
Thus the theorem applies in particular to every proper 2-edge cut of a
4-regular graph. The parity condition is needed: two triangles can be
the closures of a proper $2$-edge cut in a $6$-cycle, which has an ecd.
The graph obtained by closing the 11-vertex 2-pole in Diko's analysis
through order~18 has a proper 2-edge cut with closures $K_5$ and $K_6$
minus a perfect matching. The theorem reduces the obstruction to
$K_5$~\cite[Chapter~3, Section~3.3, discussion of the 11-vertex
multipole]{diko2026thesis}.

\begin{lemma}
\label{lem:regular-edge-replacement}
Let $k\geq2$ and let $G$ be a finite nonempty $2k$-regular multigraph.
Replace each edge $xy$ by a fresh copy of $K_{2k,2k}-ab$ and the edges
$xa$ and $yb$. The resulting graph $G^*$ is simple and $2k$-regular,
with the same edge parity and ecd index as $G$.
\end{lemma}

\begin{proof}
Each replacement preserves regularity and adds $4k^2$ edges. The final
graph is simple because the copies are disjoint and each has attachments
at distinct vertices. At each step, the two
attaching edges form a proper 2-edge cut whose closures are the graph
before replacement and $K_{2k,2k}$. The latter has an even number of
edges and ecd index $k$ by Lemma~\ref{lem:bipartite-ecd-index}. The
current graph has ecd index at least $k$ by its degree, so
Theorem~\ref{thm:two-edge-cut-decomposition} gives the same ecd index
after each replacement.
\end{proof}

Since the replacement also preserves 2-vertex-connectivity,
Questions~\ref{qu:two-connected-no-ecd} and~\ref{qu:simple-index-above-k-plus-one}
have the same answers for simple graphs and multigraphs.
Applied to the following three-vertex multigraphs, the replacement gives
simple regular graphs whose ecd index exceeds half the degree by an
arbitrarily large amount.

\begin{lemma}
\label{thm:three-vertex-multigraphs}
For every $t\geq1$, let $G_t$ be the multigraph with three vertices and
$2t$ parallel edges joining each pair of vertices. Then $G_t$ is
$4t$-regular and
\[
    \ecd(G_t)=3t=\frac{3}{4}\Delta(G_t).
\]
\end{lemma}

\begin{proof}
Every cycle of $G_t$ has length 2 or 3, so every ecd consists of $3t$
2-cycles. Any two of its cycles share a vertex and must therefore lie in
distinct colour classes. Conversely, pairing the parallel edges between each
pair of vertices gives an ecd with $3t$ colour classes. Thus
$\ecd(G_t)=3t$.
\end{proof}

\begin{theorem}
\label{thm:simple-regular-large-index}
For every integer $t\geq1$, there is a simple $4t$-regular graph $G_t^*$
with vertex connectivity 2, edge connectivity 2, and $\ecd(G_t^*)=3t$.
\end{theorem}
\begin{proof}
Apply Lemma~\ref{lem:regular-edge-replacement} with $k=2t$ to the
multigraph $G_t$ of Lemma~\ref{thm:three-vertex-multigraphs}. The resulting
graph $G_t^*$ is simple and $4t$-regular, with ecd index $3t$. Since
$G_t$ is 2-vertex-connected, so is $G_t^*$. The two attaching edges
of any copy form a 2-edge cut, and deleting its attachment vertices
$a,b$ separates the remaining vertices of the copy from the rest.
Thus both connectivities are exactly 2.
\end{proof}

\section{3-vertex-connected 4-regular graphs without an ecd}
\label{sec:three-connected-no-ecd}

We use the 4-pole $\widetilde K_5$, obtained by deleting a vertex of $K_5$
and retaining its four incident edges as dangling edges, to explain in a
unified fashion both the published constructions and the
3-vertex-connected graphs without an ecd in our census.

\paragraph{Odd passages.}
A \emph{passage} consists of two dangling edges and a path (possibly of
length zero) joining their end-vertices.
A 4-pole has the \emph{odd-passage property} if every partition of its
internal and dangling edges into even cycles and two passages uses an odd
number of internal edges in each passage. The two paths may share vertices.
Joining the four dangling edges to a new vertex gives the \emph{closure}
of the 4-pole.

Let $H$ be a 4-regular multigraph, and let $w\in V(H)$ have four
distinct neighbours. Write $H-w$ for the 4-pole obtained by deleting $w$
and retaining the four dangling edges. Then
\begin{equation}
    H-w\text{ has the odd-passage property}
    \quad\Longleftrightarrow\quad
    H\text{ has no ecd}.
    \label{eq:odd-passage-closure}
\end{equation}
Indeed, the internal graph of $H-w$ has $2|V(H-w)|-2$ edges, so the two
paths in any partition from the definition have the same parity. If both
are even, closing them through $w$ gives an ecd of $H$. Conversely, an ecd
of $H$ restricts to even cycles and two even paths in $H-w$.

The graph $K_5$ has no ecd because its 10 edges cannot be partitioned into
4-cycles, its only even cycles. Hence $\widetilde K_5$ has the odd-passage
property by~\eqref{eq:odd-passage-closure}.
\begin{proposition}[Vertex substitution]
    \label{prop:odd-pole-substitution}
    Let $B$ be a 4-regular multigraph and let
    $u\in V(B)$. For every $v\in V(B)\setminus\{u\}$, let $H_v$ be a
    4-regular multigraph without an ecd, and choose
    $w_v\in V(H_v)$ with four distinct neighbours. Replace $v$ by the
    4-pole $H_v-w_v$, attaching the four edges formerly incident with $v$
    bijectively to its dangling edges. The resulting graph has no ecd.
\end{proposition}
\begin{proof}
    Suppose that the resulting graph has an ecd. In each inserted 4-pole,
    the four dangling edges are paired by the restrictions of the cycles
    that cross the pole. This partitions the pole into two passages and
    even cycles. By~\eqref{eq:odd-passage-closure}, each passage uses an
    odd number of internal edges.

    Take a cycle $C$ of the ecd through $u$. Between successive original
    edges of $B$ on $C$ there is one internal path through an inserted
    pole, except where two original edges meet at $u$. If $C$ uses $t$
    such paths, it uses $t+1$ original edges. Each path has odd length, so
    $\lvert E(C)\rvert\equiv t+(t+1)\equiv1\pmod 2$, contradicting the
    evenness of $C$.
\end{proof}

M\'a\v{c}ajov\'a and Maz\'ak~\cite{mavcajova2013even} start with a doubled
cycle and replace every vertex except one by $\widetilde K_5$.
Liu, Cui, and Wang~\cite{liu20224} start with $K_{4,4}$ and replace every vertex
except one by a 4-pole $F_n$ formed from a chain of copies of $\widetilde K_5$.
In any partition of $F_n$ into even cycles and two passages, the internal
paths through each copy are odd by the odd-passage property of $\widetilde K_5$.
A passage made up of $t$ such paths uses $t-1$ edges between copies, so its
internal length is odd. Thus $F_n$ has the odd-passage property and its
closure has no ecd by~\eqref{eq:odd-passage-closure}.
Proposition~\ref{prop:odd-pole-substitution} applies to both
constructions. Each leaves one vertex unchanged and replaces the others
by poles of order divisible by four, so the resulting orders are all
congruent to 1 modulo 4.

The unique 4-vertex-connected graph of order 17 without an ecd in
Table~\ref{tab:ecd-4-regular}(b) is obtained from $K_5$ by replacing four
vertices by $\widetilde K_5$, leaving the fifth unchanged.

\paragraph{Paired reductions.}
A copy of $\widetilde K_5$ in $G$ consists of four vertices inducing a $K_4$,
each with one edge leaving the clique. Number these edges from 1 to 4 and
let $v_i$ be the endpoint of edge $i$ outside the clique. The endpoints may
coincide. For $i\in\{2,3,4\}$, the \emph{paired reduction} $G_{1i}$ deletes
the four internal vertices and adds two edges, one joining $v_1$ to $v_i$
and one joining the endpoints of the other pair. The degrees of the
remaining vertices are preserved. A reduction containing a loop has ecd
index $\infty$, and the ecd index of a disconnected reduction is the
maximum of its component indices.
\begin{proposition}
    \label{prop:paired-reduction}
    Let $G$ be a multigraph containing a copy of $\widetilde K_5$.
    For every $k\geq2$, $G$ admits a $k$-ecd if and only if at least one
    of its three paired reductions admits a $k$-ecd.
\end{proposition}
\begin{proof}
    A $k$-ecd of $G$ restricts to a partition of $\widetilde K_5$ into even
    cycles and two passages. Each passage has odd length, including its
    two dangling edges, by the odd-passage property. Its external endpoints
    are distinct, since otherwise it would itself be an odd cycle of the
    ecd. Replace the two passages by the added edges for their pairing
    and discard the cycles entirely within the pole. Each affected cycle
    is shortened by an even number of edges and can keep its colour.

    Conversely, relabel so that $G_{12}$ has a $k$-ecd. Replace its two
    added edges as shown below, giving each path the colour of the cycle
    that used its edge.
    
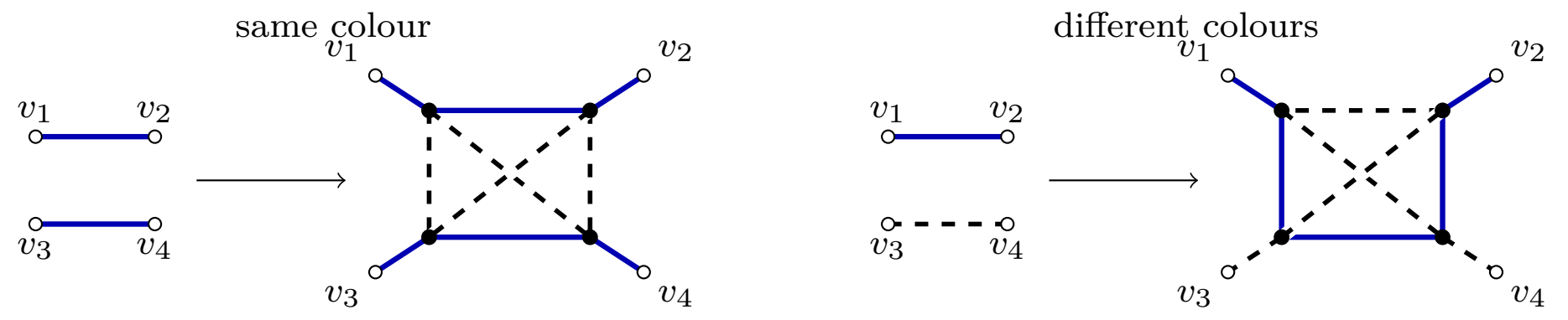
\[
    \begin{tikzpicture}[x=0.94cm,y=0.69cm,baseline=(current bounding box.center),
        same/.style={draw=blue!70!black,line width=1.2pt},
        other/.style={draw=black,dashed,line width=1.05pt},
        inner/.style={circle,fill=black,inner sep=1.25pt},
        end/.style={circle,draw=black,fill=white,inner sep=1pt},
        vlabel/.style={font=\footnotesize,inner sep=3pt},
        title/.style={font=\footnotesize}]
      \begin{scope}[shift={(0,0)}]
        \node[title] at (-0.8,1.78) {same colour};
        \draw[same] (-3.3,0.5)--(-2.3,0.5) (-3.3,-0.5)--(-2.3,-0.5);
        \foreach \x/\y/\label in {-3.3/0.5/1,-2.3/0.5/2,-3.3/-0.5/3,-2.3/-0.5/4} {
          \node[end] at (\x,\y) {};
          \ifdim\y pt>0pt
            \node[vlabel,above] at (\x,\y) {$v_{\label}$};
          \else
            \node[vlabel,below] at (\x,\y) {$v_{\label}$};
          \fi
        }
        \draw[->] (-1.95,0)--(-0.7,0);
        \coordinate (x1) at (0,0.8);
        \coordinate (x2) at (1.35,0.8);
        \coordinate (x3) at (0,-0.65);
        \coordinate (x4) at (1.35,-0.65);
        \draw[other] (x1)--(x3) (x2)--(x4);
        \draw[other] (x1)--(x4);
        \draw[other,preaction={draw=white,line width=3pt}] (x2)--(x3);
        \draw[same] (x1)--(x2) (x3)--(x4);
        \draw[same] (x1)--(-0.45,1.2) (x2)--(1.8,1.2)
                    (x3)--(-0.45,-1.05) (x4)--(1.8,-1.05);
        \foreach \x in {x1,x2,x3,x4} {
          \node[inner] at (\x) {};
        }
        \node[end] at (-0.45,1.2) {};
        \node[end] at (1.8,1.2) {};
        \node[end] at (-0.45,-1.05) {};
        \node[end] at (1.8,-1.05) {};
        \node[vlabel,above left] at (-0.45,1.2) {$v_1$};
        \node[vlabel,above right] at (1.8,1.2) {$v_2$};
        \node[vlabel,below left] at (-0.45,-1.05) {$v_3$};
        \node[vlabel,below right] at (1.8,-1.05) {$v_4$};
      \end{scope}
      \begin{scope}[shift={(7.15,0)}]
        \node[title] at (-0.8,1.78) {different colours};
        \draw[same] (-3.3,0.5)--(-2.3,0.5);
        \draw[other] (-3.3,-0.5)--(-2.3,-0.5);
        \foreach \x/\y/\label in {-3.3/0.5/1,-2.3/0.5/2,-3.3/-0.5/3,-2.3/-0.5/4} {
          \node[end] at (\x,\y) {};
          \ifdim\y pt>0pt
            \node[vlabel,above] at (\x,\y) {$v_{\label}$};
          \else
            \node[vlabel,below] at (\x,\y) {$v_{\label}$};
          \fi
        }
        \draw[->] (-1.95,0)--(-0.7,0);
        \coordinate (x1) at (0,0.8);
        \coordinate (x2) at (1.35,0.8);
        \coordinate (x3) at (0,-0.65);
        \coordinate (x4) at (1.35,-0.65);
        \draw[same] (x1)--(x3) (x3)--(x4) (x4)--(x2);
        \draw[other] (x2)--(x1);
        \draw[other] (x1)--(x4);
        \draw[other,preaction={draw=white,line width=3pt}] (x3)--(x2);
        \draw[same] (x1)--(-0.45,1.2) (x2)--(1.8,1.2);
        \draw[other] (x3)--(-0.45,-1.05) (x4)--(1.8,-1.05);
        \foreach \x in {x1,x2,x3,x4} {
          \node[inner] at (\x) {};
        }
        \node[end] at (-0.45,1.2) {};
        \node[end] at (1.8,1.2) {};
        \node[end] at (-0.45,-1.05) {};
        \node[end] at (1.8,-1.05) {};
        \node[vlabel,above left] at (-0.45,1.2) {$v_1$};
        \node[vlabel,above right] at (1.8,1.2) {$v_2$};
        \node[vlabel,below left] at (-0.45,-1.05) {$v_3$};
        \node[vlabel,below right] at (1.8,-1.05) {$v_4$};
      \end{scope}
    \end{tikzpicture}
    \]
    The paths have odd length and share internal vertices only when their
    colours differ. Give the unused $4$-cycle in the left panel a second
    colour. The result is a $k$-ecd of $G$.
\end{proof}

In particular, if $G$ is 4-regular, its ecd index is the smallest of the
three paired-reduction indices.

\paragraph{Census obstructions.}
Each graph without an ecd in Tables~\ref{tab:ecd-4-regular}
and~\ref{tab:ecd-4-regular-multigraphs} has a reduction tree with steps at
$2$-edge cuts (Theorem~\ref{thm:two-edge-cut-decomposition}) or at copies of
$\widetilde K_5$ (Proposition~\ref{prop:paired-reduction}).
Its leaves are $K_5$ or graphs with a loop or a cut vertex
(Proposition~\ref{prop:4reg-cut-vertex-no-ecd}). Some trees use both types
of reduction. The same obstructions explain the 144 order-20 cases left
open by the reductions used in the
thesis~\cite[Chapter~3, subsection ``3-connected 4-regular graphs'']{diko2026thesis}.

Up to order 20, the census contains 3-vertex-connected 4-regular graphs
without an ecd only at orders congruent to 1 modulo 4. Joining two to six
4-poles obtained from such graphs through order 17 also produced no
examples at other orders~\cite[Chapter~3, subsection ``3-connected 4-regular graphs'']{diko2026thesis}.
\begin{question}
    \label{qu:three-connected-mod-4}
    Does there exist a 3-vertex-connected 4-regular graph without an ecd where
    $\lvert G \rvert \not\equiv 1 \pmod 4$?
\end{question}

Liu, Cui, and Wang~\cite{liu20224} constructed infinitely many
2-vertex-connected 4-regular graphs of girth at least four without an ecd.
Every 3-vertex-connected simple 4-regular graph without an ecd in our census
contains a $K_4$. Thus any example in the following question of
Markstr\"om~\cite{MARKSTROM20122676} must have order greater than 20.
\begin{question}
    \label{qu:three-connected-triangle-free}
    Does there exist a 3-vertex-connected 4-regular graph of girth at least
    four without an ecd?
\end{question}

\section{Circulants}
\label{sec:circulants}
A \emph{circulant graph} $C_n(j_1,\ldots,j_k)$ has vertex set
$\mathbb Z_n$, with $v$ adjacent to $v\pm j_i$ for every $i$. We assume
$0<j_1<\cdots<j_k<n/2$, so the graph is simple and $2k$-regular. It is
connected if and only if $\gcd(n,j_1,\ldots,j_k)=1$
\cite{heuberger2003planarity}.

\begin{proposition}
\label{prop:circulant-parity}
Let $G$ be a simple circulant.
\begin{center}
\begin{tabular}{@{}l@{\enspace}lll@{}}
\hline
\multicolumn{2}{@{}l}{\textbf{Degree of $G$}} & \textbf{Even order and connected} & \textbf{Odd order} \\
\hline
$4r$, & $r\geq1$ & $\ecd(G)=2r$ & $\ecd(G)\geq 2r+1$ if $G$ has an ecd \\
$4r+2$, & $r\geq0$ & $\ecd(G)=2r+1$ & $G$ has no ecd \\
\hline
\end{tabular}
\end{center}
\end{proposition}

\begin{proof}
Write $n=|V(G)|$. If $n$ is even and $G$ is connected, then $G$ is a connected Cayley graph of the
abelian group $\mathbb Z_n$ and is therefore 1-factorizable by Stong's
theorem \cite{stong19851}.
Proposition~\ref{prop:ecd-index-chromatic-index} then gives $\ecd(G)$ equal
to half the degree. If $n$ is odd and the degree is $4r+2$, then
$|E|=(2r+1)n$ is odd, whereas every ecd has an even number of edges.
If $n$ is odd and the degree is $4r$, Proposition~\ref{prop:ecd-index-chromatic-index}
gives $\ecd(G)\geq 2r+1$.
\end{proof}

Thus the only unresolved parity case has odd order and degree divisible by
4. We conjecture the following.
\begin{conjecture}
For every $r\geq1$, every connected simple $4r$-regular circulant of odd
order has ecd index $2r+1$, with the exception of $K_5$, which has no ecd.
\end{conjecture}
The rest of this section gathers evidence for this conjecture.
Theorem~\ref{thm:circulant-two-hamiltonian-cycles} resolves it for many
4-regular circulants, and Theorem~\ref{thm:circulant-12-13} for the
two-step families $C_n(1,2)$ and $C_n(1,3)$.

\begin{theorem}
\label{thm:circulant-two-hamiltonian-cycles}
Let $C_n(a,b)$ be a simple 4-regular circulant of odd order $n$
with $a$ and $b$ invertible modulo $n$. Let $s,t\in\{1,\ldots,n-1\}$ be
the unique integers such that
\[
    s\equiv ba^{-1}\pmod n
    \qquad\text{and}\qquad
    t\equiv ab^{-1}\pmod n.
\]
If $s$ and $t$ have the same parity, then
\[
    \ecd(C_n(a,b))=3.
\]
Moreover, the 4-cycle $(0,a,a+b,b)$ can be one of the three colour classes.
\end{theorem}

\begin{proof}
Let $Q=(0,a,a+b,b)$ and remove its four edges. Since $a$ is invertible
modulo $n$, the edges of step $a$ form a Hamilton cycle. Removing the two
edges of $Q$
splits it into paths with common endpoint pairs $\{a,b\}$ and
$\{a+b,0\}$ and lengths $s-1$ and $n-s-1$. The step-$b$ Hamilton cycle
similarly splits into paths with the same endpoint pairs and lengths $t-1$
and $n-t-1$. Indeed, starting at $b$, after $t-1$ further step-$b$ edges
we reach $tb=a$, while starting at $a+b$, after $n-t-1$ such edges we
reach $a+(n-t)b=0$, with all equalities taken modulo $n$.

Colour each of the four paths alternately red and blue. For each pair of
paths with the same endpoints, choose the colours at one endpoint to be equal.
Since the two path lengths have the same parity, their colours at the other
endpoint are also equal. Hence every vertex has degree zero or two in each
colour, and the two incident edges of one colour have different steps. The red
and blue edges therefore form cycles alternating between steps $a$ and $b$;
these cycles are even. Together with $Q$, they form a $3$-ecd.
Proposition~\ref{prop:ecd-index-chromatic-index} gives index at least 3,
so the index is exactly 3.
\end{proof}

The families $C_n(1,2)$ and $C_n(1,3)$ are not covered for every order by
Theorem~\ref{thm:circulant-two-hamiltonian-cycles}: for instance,
$C_n(1,2)$ with $n\equiv1\pmod4$ and $C_n(1,3)$ with $n\equiv5\pmod6$
fall outside its hypotheses. They nevertheless admit particularly simple
inductive constructions.

\begin{theorem}
\label{thm:circulant-12-13}
For every odd $n\geq7$,
\[
    \ecd(C_n(1,2))=\ecd(C_n(1,3))=3.
\]
\end{theorem}

\begin{figure}[p]
    \centering

    \centering
    \begin{tikzpicture}[scale=0.82]
        \foreach \i in {0,...,5} {
            \node[draw, circle, fill=black, inner sep=0pt, minimum size=0.15cm, label=below:$\i$] (v\i) at (\i, 0) {};
        }

        \draw[myedge, red] (v0) -- node[midway, above, text=red, fill=white, inner sep=1pt, font=\footnotesize, outer sep=1pt] {5} (v1);
        \draw[myedge, green!50!black] (v1) -- node[midway, above, text=green!50!black, fill=white, inner sep=1pt, font=\footnotesize, outer sep=1pt] {3} (v2);
        \draw[myedge, red] (v2) -- node[midway, above, text=red, fill=white, inner sep=1pt, font=\footnotesize, outer sep=1pt] {5} (v3);
        \draw[myedge, green!50!black] (v3) -- node[midway, above, text=green!50!black, fill=white, inner sep=1pt, font=\footnotesize, outer sep=1pt] {3} (v4);
        \draw[myedge, red] (v4) -- node[midway, above, text=red, fill=white, inner sep=1pt, font=\footnotesize, outer sep=1pt] {5} (v5);

        \draw[myedge, red] (v0) to[out=80, in=100] node[pos=0.5, above, text=red, font=\footnotesize, fill=white, inner sep=1pt] {6} (v3);

        \draw[myedge, green!50!black] (v1) to[out=80, in=100] node[pos=0.5, above, text=green!50!black, font=\footnotesize, fill=white, inner sep=1pt] {4} (v4);

        \draw[myedge, green!50!black] (v2) to[out=80, in=100] node[pos=0.5, above, text=green!50!black, font=\footnotesize, fill=white, inner sep=1pt] {4} (v5);


        \draw[myedge, green!50!black] (v0) -- node[midway, above, text=green!50!black, font=\footnotesize] {3} ++(-0.8, 0) node[left, text=black, font=\footnotesize] {$e_1$};
        \draw[myedge, green!50!black] (v0) to[out=120, in=0] node[midway, above, text=green!50!black, font=\footnotesize] {4} ++(-1.0, 0.5) node[left, text=black, font=\footnotesize] {$e_2$};
        \draw[myedge, red] (v1) to[out=120, in=0] node[midway, above, text=red, font=\footnotesize] {6} ++(-2, 0.8) node[left, text=black, font=\footnotesize] {$e_3$};
        \draw[myedge, red] (v2) to[out=120, in=0] node[midway, above, text=red, font=\footnotesize] {6} ++(-3, 1.3) node[left, text=black, font=\footnotesize] {$e_4$};

        \draw[myedge, green!50!black] (v5) -- node[midway, above, text=green!50!black, font=\footnotesize] {3} ++(0.8, 0) node[right, text=black, font=\footnotesize] {$e_1$};
        \draw[myedge, green!50!black] (v3) to[out=60, in=180] node[midway, above, text=green!50!black, font=\footnotesize] {4} ++(3, 1.3) node[right, text=black, font=\footnotesize] {$e_2$};
        \draw[myedge, red] (v4) to[out=60, in=180] node[midway, above, text=red, font=\footnotesize] {6} ++(2, 0.8) node[right, text=black, font=\footnotesize] {$e_3$};
        \draw[myedge, red] (v5) to[out=60, in=180] node[midway, above, text=red, font=\footnotesize] {6} ++(1.0, 0.5) node[right, text=black, font=\footnotesize] {$e_4$};

    \end{tikzpicture}
    \par\smallskip
    {\small (a) Six-vertex extension strip.}
    \par\medskip
    \begin{tikzpicture}[scale=0.82]
        \foreach \i in {0,...,6} {
            \node[draw, circle, fill=black, inner sep=0pt, minimum size=0.15cm, label=below:$\i$] (v\i) at (\i, 0) {};
        }

        \draw[myedge, blue] (v0) -- node[mylabel] {1} (v1);
        \draw[myedge, red] (v1) -- node[midway, above, text=red, fill=white, inner sep=1pt, font=\footnotesize, outer sep=1pt] {5} (v2);
        \draw[myedge, green!50!black] (v2) -- node[midway, above, text=green!50!black, fill=white, inner sep=1pt, font=\footnotesize, outer sep=1pt] {3} (v3);
        \draw[myedge, blue] (v3) -- node[mylabel] {1} (v4);
        \draw[myedge, green!50!black] (v4) -- node[midway, above, text=green!50!black, fill=white, inner sep=1pt, font=\footnotesize, outer sep=1pt] {3} (v5);
        \draw[myedge, red] (v5) -- node[midway, above, text=red, fill=white, inner sep=1pt, font=\footnotesize, outer sep=1pt] {5} (v6);

        \draw[myedge, blue] (v0) to[out=70, in=110, looseness=1.2] node[mylabel] {2} (v3);
        \draw[myedge, blue] (v1) to[out=70, in=110, looseness=1.2] node[mylabel] {2} (v4);
        \draw[myedge, green!50!black] (v2) to[out=70, in=110, looseness=1.2] node[midway, above, text=green!50!black, fill=white, inner sep=1pt, font=\footnotesize, outer sep=1pt] {4} (v5);
        \draw[myedge, green!50!black] (v3) to[out=70, in=110, looseness=1.2] node[midway, above, text=green!50!black, fill=white, inner sep=1pt, font=\footnotesize, outer sep=1pt] {4} (v6);

        \draw[myedge, green!50!black] (v0) -- node[midway, above, text=green!50!black, font=\footnotesize] {3} ++(-0.8, 0) node[left, text=black, font=\footnotesize] {$e_1$};
        \draw[myedge, green!50!black] (v0) to[out=120, in=0] node[midway, above, text=green!50!black, font=\footnotesize] {4} ++(-1.0, 0.5) node[left, text=black, font=\footnotesize] {$e_2$};
        \draw[myedge, red] (v1) to[out=120, in=0] node[midway, above, text=red, font=\footnotesize] {6} ++(-2, 0.8) node[left, text=black, font=\footnotesize] {$e_3$};
        \draw[myedge, red] (v2) to[out=120, in=0] node[midway, above, text=red, font=\footnotesize] {6} ++(-3,1.3) node[left, text=black, font=\footnotesize] {$e_4$};

        \draw[myedge, green!50!black] (v6) -- node[midway, above, text=green!50!black, font=\footnotesize] {3} ++(0.8, 0) node[right, text=black, font=\footnotesize] {$e_1$};
        \draw[myedge, green!50!black] (v4) to[out=60, in=180] node[midway, above, text=green!50!black, font=\footnotesize] {4} ++(3, 1.3) node[right, text=black, font=\footnotesize] {$e_2$};
        \draw[myedge, red] (v5) to[out=60, in=180] node[midway, above, text=red, font=\footnotesize] {6} ++(2, 0.8) node[right, text=black, font=\footnotesize] {$e_3$};
        \draw[myedge, red] (v6) to[out=60, in=180] node[midway, above, text=red, font=\footnotesize] {6} ++(1.0, 0.5) node[right, text=black, font=\footnotesize] {$e_4$};

    \end{tikzpicture}
    \par\smallskip
    {\small (b) Base strip of order 7.}
    \par\medskip

    \centering
    \begin{tikzpicture}[scale=0.82]
        \foreach \i in {0,...,8} {
            \node[draw, circle, fill=black, inner sep=0pt, minimum size=0.15cm, label=below:$\i$] (v\i) at (\i, 0) {};
        }

        \draw[myedge, red] (v0) -- node[midway, above, text=red, fill=white, inner sep=1pt, font=\footnotesize, outer sep=1pt] {5} (v1);
        \draw[myedge, blue] (v1) -- node[mylabel] {1} (v2);
        \draw[myedge, red] (v2) -- node[midway, above, text=red, fill=white, inner sep=1pt, font=\footnotesize, outer sep=1pt] {5} (v3);
        \draw[myedge, blue] (v3) -- node[mylabel] {1} (v4);
        \draw[myedge, green!50!black] (v4) -- node[midway, above, text=green!50!black, fill=white, inner sep=1pt, font=\footnotesize, outer sep=1pt] {3} (v5);
        \draw[myedge, blue] (v5) -- node[mylabel] {1} (v6);
        \draw[myedge, green!50!black] (v6) -- node[midway, above, text=green!50!black, fill=white, inner sep=1pt, font=\footnotesize, outer sep=1pt] {3} (v7);
        \draw[myedge, red] (v7) -- node[midway, above, text=red, fill=white, inner sep=1pt, font=\footnotesize, outer sep=1pt] {5} (v8);

        \draw[myedge, red] (v0) to[out=70, in=110, looseness=1.2] node[midway, above, text=red, fill=white, inner sep=1pt, font=\footnotesize, outer sep=1pt] {6} (v3);
        \draw[myedge, blue] (v1) to[out=70, in=110, looseness=1.2] node[mylabel] {2} (v4);
        \draw[myedge, blue] (v2) to[out=70, in=110, looseness=1.2] node[mylabel] {2} (v5);
        \draw[myedge, blue] (v3) to[out=70, in=110, looseness=1.2] node[mylabel] {2} (v6);
        \draw[myedge, green!50!black] (v4) to[out=70, in=110, looseness=1.2] node[midway, above, text=green!50!black, fill=white, inner sep=1pt, font=\footnotesize, outer sep=1pt] {4} (v7);
        \draw[myedge, green!50!black] (v5) to[out=70, in=110, looseness=1.2] node[midway, above, text=green!50!black, fill=white, inner sep=1pt, font=\footnotesize, outer sep=1pt] {4} (v8);

        \draw[myedge, green!50!black] (v0) -- node[midway, above, text=green!50!black, font=\footnotesize] {3} ++(-0.8, 0) node[left, text=black, font=\footnotesize] {$e_1$};
        \draw[myedge, green!50!black] (v0) to[out=120, in=0] node[midway, above, text=green!50!black, font=\footnotesize] {4} ++(-1.0, 0.5) node[left, text=black, font=\footnotesize] {$e_2$};
        \draw[myedge, red] (v1) to[out=120, in=0] node[midway, above, text=red, font=\footnotesize] {6} ++(-2, 0.8) node[left, text=black, font=\footnotesize] {$e_3$};
        \draw[myedge, red] (v2) to[out=120, in=0] node[midway, above, text=red, font=\footnotesize] {6} ++(-3, 1.3) node[left, text=black, font=\footnotesize] {$e_4$};

        \draw[myedge, green!50!black] (v8) -- node[midway, above, text=green!50!black, font=\footnotesize] {3} ++(0.8, 0) node[right, text=black, font=\footnotesize] {$e_1$};
        \draw[myedge, green!50!black] (v6) to[out=60, in=180] node[midway, above, text=green!50!black, font=\footnotesize] {4} ++(3, 1.3) node[right, text=black, font=\footnotesize] {$e_2$};
        \draw[myedge, red] (v7) to[out=60, in=180] node[midway, above, text=red, font=\footnotesize] {6} ++(2, 0.8) node[right, text=black, font=\footnotesize] {$e_3$};
        \draw[myedge, red] (v8) to[out=60, in=180] node[midway, above, text=red, font=\footnotesize] {6} ++(1.0, 0.5) node[right, text=black, font=\footnotesize] {$e_4$};

    \end{tikzpicture}
    \par\smallskip
    {\small (c) Base strip of order 9.}
    \par\medskip

    \centering
    \begin{tikzpicture}[scale=0.82]
        \foreach \i in {0,...,10} {
            \node[draw, circle, fill=black, inner sep=0pt, minimum size=0.15cm, label=below:$\i$] (v\i) at (\i, 0) {};
        }

        \draw[myedge, red] (v0) -- node[midway, above, text=red, fill=white, inner sep=1pt, font=\footnotesize, outer sep=1pt] {5} (v1);
        \draw[myedge, blue] (v1) -- node[mylabel] {1} (v2);
        \draw[myedge, red] (v2) -- node[midway, above, text=red, fill=white, inner sep=1pt, font=\footnotesize, outer sep=1pt] {5} (v3);
        \draw[myedge, blue] (v3) -- node[mylabel] {1} (v4);
        \draw[myedge, red] (v4) -- node[midway, above, text=red, fill=white, inner sep=1pt, font=\footnotesize, outer sep=1pt] {5} (v5);
        \draw[myedge, blue] (v5) -- node[mylabel] {1} (v6);
        \draw[myedge, green!50!black] (v6) -- node[midway, above, text=green!50!black, fill=white, inner sep=1pt, font=\footnotesize, outer sep=1pt] {3} (v7);
        \draw[myedge, red] (v7) -- node[midway, above, text=red, fill=white, inner sep=1pt, font=\footnotesize, outer sep=1pt] {5} (v8);
        \draw[myedge, green!50!black] (v8) -- node[midway, above, text=green!50!black, fill=white, inner sep=1pt, font=\footnotesize, outer sep=1pt] {3} (v9);
        \draw[myedge, red] (v9) -- node[midway, above, text=red, fill=white, inner sep=1pt, font=\footnotesize, outer sep=1pt] {5} (v10);

        \draw[myedge, red] (v0) to[out=70, in=110, looseness=1.2] node[midway, above, text=red, fill=white, inner sep=1pt, font=\footnotesize, outer sep=1pt] {6} (v3);
        \draw[myedge, blue] (v1) to[out=70, in=110, looseness=1.2] node[mylabel] {2} (v4);
        \draw[myedge, blue] (v2) to[out=70, in=110, looseness=1.2] node[mylabel] {2} (v5);
        \draw[myedge, blue] (v3) to[out=70, in=110, looseness=1.2] node[mylabel] {2} (v6);
        \draw[myedge, red] (v4) to[out=70, in=110, looseness=1.2] node[midway, above, text=red, fill=white, inner sep=1pt, font=\footnotesize, outer sep=1pt] {6} (v7);
        \draw[myedge, red] (v5) to[out=70, in=110, looseness=1.2] node[midway, above, text=red, fill=white, inner sep=1pt, font=\footnotesize, outer sep=1pt] {6} (v8);
        \draw[myedge, green!50!black] (v6) to[out=70, in=110, looseness=1.2] node[midway, above, text=green!50!black, fill=white, inner sep=1pt, font=\footnotesize, outer sep=1pt] {4} (v9);
        \draw[myedge, green!50!black] (v7) to[out=70, in=110, looseness=1.2] node[midway, above, text=green!50!black, fill=white, inner sep=1pt, font=\footnotesize, outer sep=1pt] {4} (v10);

        \draw[myedge, green!50!black] (v0) -- node[midway, above, text=green!50!black, font=\footnotesize] {3} ++(-0.8, 0) node[left, text=black, font=\footnotesize] {$e_1$};
        \draw[myedge, green!50!black] (v0) to[out=120, in=0] node[midway, above, text=green!50!black, font=\footnotesize] {4} ++(-1.0, 0.5) node[left, text=black, font=\footnotesize] {$e_2$};
        \draw[myedge, red] (v1) to[out=120, in=0] node[midway, above, text=red, font=\footnotesize] {6} ++(-2, 0.8) node[left, text=black, font=\footnotesize] {$e_3$};
        \draw[myedge, red] (v2) to[out=120, in=0] node[midway, above, text=red, font=\footnotesize] {6} ++(-3, 1.3) node[left, text=black, font=\footnotesize] {$e_4$};

        \draw[myedge, green!50!black] (v10) -- node[midway, above, text=green!50!black, font=\footnotesize] {3} ++(0.8, 0) node[right, text=black, font=\footnotesize] {$e_1$};
        \draw[myedge, green!50!black] (v8) to[out=60, in=180] node[midway, above, text=green!50!black, font=\footnotesize] {4} ++(3, 1.3) node[right, text=black, font=\footnotesize] {$e_2$};
        \draw[myedge, red] (v9) to[out=60, in=180] node[midway, above, text=red, font=\footnotesize] {6} ++(2, 0.8) node[right, text=black, font=\footnotesize] {$e_3$};
        \draw[myedge, red] (v10) to[out=60, in=180] node[midway, above, text=red, font=\footnotesize] {6} ++(1.0, 0.5) node[right, text=black, font=\footnotesize] {$e_4$};

    \end{tikzpicture}
    \par\smallskip
    {\small (d) Base strip of order 11.}
    \caption{Base and extension strips for the 3-ecd colouring of
    $C_n(1,3)$. Joining dangling edges with the same boundary name $e_i$
    closes a base; their numeric colour labels agree.
    Inserting panel (a) increases the order by 6.}
    \label{fig:circulant-13-construction}
\end{figure}
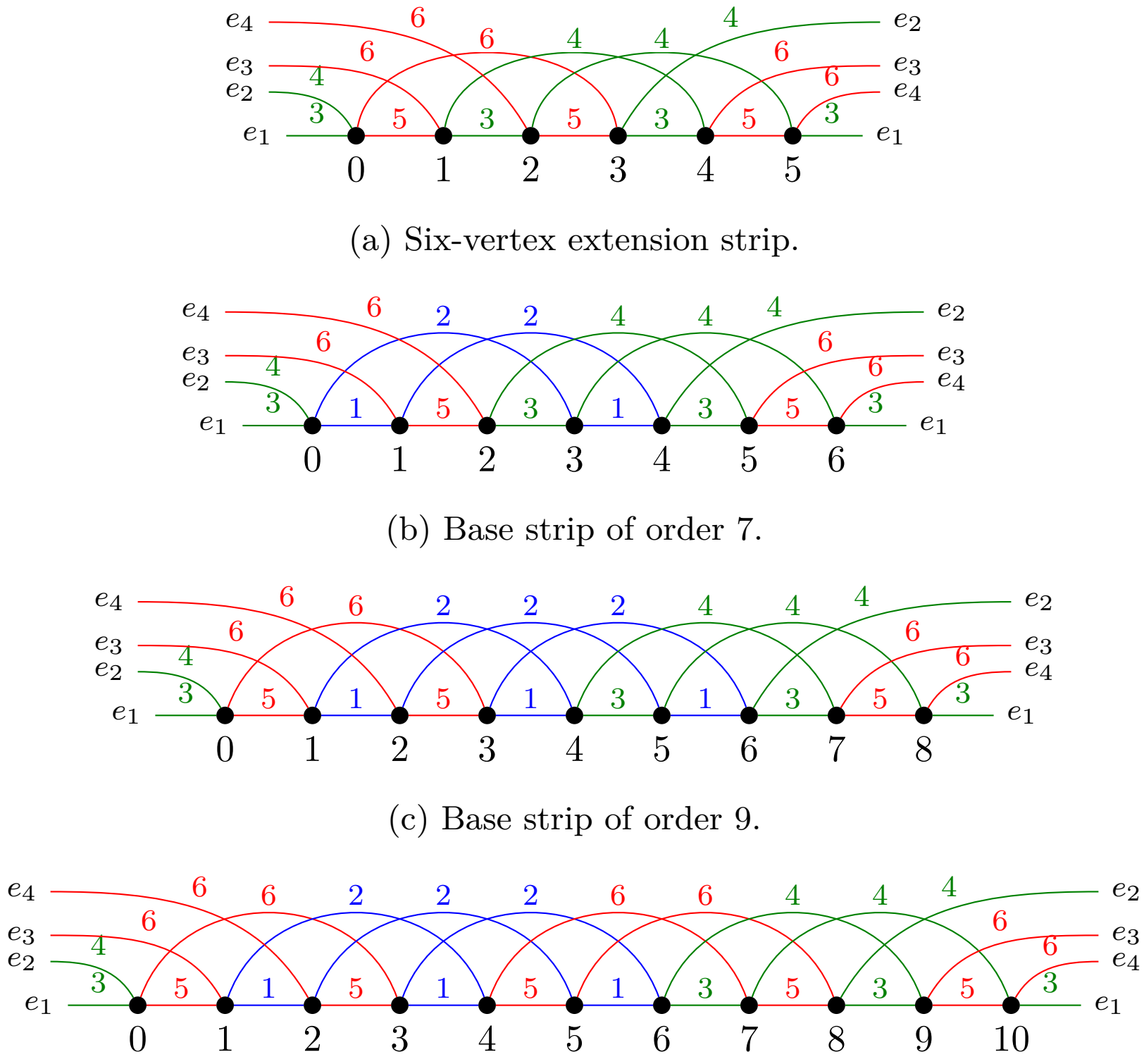

\begin{proof}
For $C_n(1,2)$, the bases of orders 7 and 9 and the four-vertex extension
in Figure~\ref{fig:circulant-12-construction} cover the two odd residue
classes modulo 4. For $C_n(1,3)$, the bases of orders 7, 9, and 11 and the
six-vertex extension in Figure~\ref{fig:circulant-13-construction} cover the
three odd residue classes modulo 6. In both figures, the edges labelled 1
and 2 form one colour class, those labelled 3 and 4 form a second, and those
labelled 5 and 6 form a third; joining dangling edges with the same boundary
name $e_i$ closes a strip, while inserting the displayed extension strip preserves its
boundary colouring. In each case the local conditions of
Proposition~\ref{prop:ecd-colouring-characterization} can be read directly from
the figure. Full details of the induction are given in the thesis
\cite[Chapter~4, Section~1 and Appendix~D]{diko2026thesis}. Proposition
\ref{prop:ecd-index-chromatic-index} gives the lower bound 3.
\end{proof}

\begin{figure}
    \centering
    \begin{tikzpicture}[scale=1.0]
        \foreach \i in {0,...,6} {
            \node[draw, circle, fill=black, inner sep=0pt, minimum size=0.15cm, label=below:$\i$] (v\i) at (\i, 0) {};
        }

        \draw[myedge, red] (v0) -- node[midway, above, text=red, fill=white, inner sep=1pt, font=\footnotesize, outer sep=1pt] {5} (v1);
        \draw[myedge, green!50!black] (v1) -- node[midway, above, text=green!50!black, fill=white, inner sep=1pt, font=\footnotesize, outer sep=1pt] {4} (v2);
        \draw[myedge, blue] (v2) -- node[mylabel] {1} (v3);
        \draw[myedge, green!50!black] (v3) -- node[midway, above, text=green!50!black, fill=white, inner sep=1pt, font=\footnotesize, outer sep=1pt] {4} (v4);
        \draw[myedge, blue] (v4) -- node[mylabel] {1} (v5);
        \draw[myedge, red] (v5) -- node[midway, above, text=red, fill=white, inner sep=1pt, font=\footnotesize, outer sep=1pt] {5} (v6);

        \draw[myedge, green!50!black] (v0) to[out=60, in=120] node[midway, above, text=green!50!black, fill=white, inner sep=1pt, font=\footnotesize, outer sep=1pt] {3} (v2);
        \draw[myedge, green!50!black] (v1) to[out=60, in=120] node[midway, above, text=green!50!black, fill=white, inner sep=1pt, font=\footnotesize, outer sep=1pt] {3} (v3);
        \draw[myedge, blue] (v2) to[out=60, in=120] node[mylabel] {2} (v4);
        \draw[myedge, blue] (v3) to[out=60, in=120] node[mylabel] {2} (v5);
        \draw[myedge, green!50!black] (v4) to[out=60, in=120] node[midway, above, text=green!50!black, fill=white, inner sep=1pt, font=\footnotesize, outer sep=1pt] {3} (v6);

        \draw[myedge, green!50!black] (v0) -- node[midway, above, text=green!50!black, font=\footnotesize] {4} ++(-0.8, 0) node[left, text=black, font=\footnotesize] {$e_1$};
        \draw[myedge, red] (v0) to[out=120, in=0] node[midway, above, text=red, font=\footnotesize] {6} ++(-1, 0.5) node[left, text=black, font=\footnotesize] {$e_2$};
        \draw[myedge, red] (v1) to[out=120, in=0] node[midway, above, text=red, font=\footnotesize] {6} ++(-1.5, 0.8) node[left, text=black, font=\footnotesize] {$e_3$};

        \draw[myedge, green!50!black] (v6) -- node[midway, above, text=green!50!black, font=\footnotesize] {4} ++(0.8, 0) node[right, text=black, font=\footnotesize] {$e_1$};
        \draw[myedge, red] (v5) to[out=60, in=180] node[midway, above, text=red, font=\footnotesize] {6} ++(1.5, 0.8) node[right, text=black, font=\footnotesize] {$e_2$};
        \draw[myedge, red] (v6) to[out=60, in=180] node[midway, above, text=red, font=\footnotesize] {6} ++(1, 0.5) node[right, text=black, font=\footnotesize] {$e_3$};

    \end{tikzpicture}
    \par\smallskip
    {\small (a) Base strip of order 7.}
    \par\medskip

    \begin{tikzpicture}[scale=1.0]
        \foreach \i in {0,...,3} {
             \node[draw, circle, fill=black, inner sep=0pt, minimum size=0.15cm, label=below:$\i$] (v\i) at (\i, 0) {};
        }

        \draw[myedge, green!50!black] (v0) -- node[midway, above, text=green!50!black, fill=white, inner sep=1pt, font=\footnotesize, outer sep=1pt] {3} (v1);
        \draw[myedge, green!50!black] (v1) -- node[midway, above, text=green!50!black, fill=white, inner sep=1pt, font=\footnotesize, outer sep=1pt] {4} (v2);
        \draw[myedge, green!50!black] (v2) -- node[midway, above, text=green!50!black, fill=white, inner sep=1pt, font=\footnotesize, outer sep=1pt] {3} (v3);

        \draw[myedge, red] (v0) to[out=80, in=100] node[pos=0.5, above, text=red, font=\footnotesize, fill=white, inner sep=1pt] {5} (v2);

        \draw[myedge, red] (v1) to[out=80, in=100] node[pos=0.5, above, text=red, font=\footnotesize, fill=white, inner sep=1pt] {5} (v3);


        \draw[myedge, green!50!black] (v0) -- node[midway, above, text=green!50!black, font=\footnotesize] {4} ++(-0.8, 0) node[left, text=black, font=\footnotesize] {$e_1$};
        \draw[myedge, red] (v0) to[out=120, in=0] node[midway, above, text=red, font=\footnotesize] {6} ++(-1.0, 0.5) node[left, text=black, font=\footnotesize] {$e_2$};
        \draw[myedge, red] (v1) to[out=120, in=0] node[midway, above, text=red, font=\footnotesize] {6} ++(-1.7, 0.8) node[left, text=black, font=\footnotesize] {$e_3$};
        \draw[myedge, green!50!black] (v3) -- node[midway, above, text=green!50!black, font=\footnotesize] {4} ++(0.8, 0) node[right, text=black, font=\footnotesize] {$e_1$};
        \draw[myedge, red] (v3) to[out=60, in=180] node[midway, above, text=red, font=\footnotesize] {6} ++(1.0, 0.5) node[right, text=black, font=\footnotesize] {$e_3$};
        \draw[myedge, red] (v2) to[out=60, in=180] node[midway, above, text=red, font=\footnotesize] {6}  ++(1.7, 0.8) node[right, text=black, font=\footnotesize] {$e_2$};
    \end{tikzpicture}
    \par\smallskip
    {\small (b) Four-vertex extension strip.}
    \par\medskip

    \begin{tikzpicture}[scale=1.0]
        \foreach \i in {0,...,8} {
            \node[draw, circle, fill=black, inner sep=0pt, minimum size=0.15cm, label=below:$\i$] (v\i) at (\i, 0) {};
        }

        \draw[myedge, green!50!black] (v0) -- node[midway, above, text=green!50!black, fill=white, inner sep=1pt, font=\footnotesize, outer sep=1pt] {3} (v1);
        \draw[myedge, green!50!black] (v1) -- node[midway, above, text=green!50!black, fill=white, inner sep=1pt, font=\footnotesize, outer sep=1pt] {4} (v2);
        \draw[myedge, red] (v2) -- node[midway, above, text=red, fill=white, inner sep=1pt, font=\footnotesize, outer sep=1pt] {6} (v3);
        \draw[myedge, green!50!black] (v3) -- node[midway, above, text=green!50!black, fill=white, inner sep=1pt, font=\footnotesize, outer sep=1pt] {4} (v4);
        \draw[myedge, blue] (v4) -- node[mylabel] {1} (v5);
        \draw[myedge, green!50!black] (v5) -- node[midway, above, text=green!50!black, fill=white, inner sep=1pt, font=\footnotesize, outer sep=1pt] {4} (v6);
        \draw[myedge, blue] (v6) -- node[mylabel] {1} (v7);
        \draw[myedge, red] (v7) -- node[midway, above, text=red, fill=white, inner sep=1pt, font=\footnotesize, outer sep=1pt] {5} (v8);

        \draw[myedge, red] (v0) to[out=60, in=120] node[midway, above, text=red, fill=white, inner sep=1pt, font=\footnotesize, outer sep=1pt] {5} (v2);
        \draw[myedge, red] (v1) to[out=60, in=120] node[midway, above, text=red, fill=white, inner sep=1pt, font=\footnotesize, outer sep=1pt] {5} (v3);
        \draw[myedge, green!50!black] (v2) to[out=60, in=120] node[midway, above, text=green!50!black, fill=white, inner sep=1pt, font=\footnotesize, outer sep=1pt] {3} (v4);
        \draw[myedge, green!50!black] (v3) to[out=60, in=120] node[midway, above, text=green!50!black, fill=white, inner sep=1pt, font=\footnotesize, outer sep=1pt] {3} (v5);
        \draw[myedge, blue] (v4) to[out=60, in=120] node[mylabel] {2} (v6);
        \draw[myedge, blue] (v5) to[out=60, in=120] node[mylabel] {2} (v7);
        \draw[myedge, green!50!black] (v6) to[out=60, in=120] node[midway, above, text=green!50!black, fill=white, inner sep=1pt, font=\footnotesize, outer sep=1pt] {3} (v8);

        \draw[myedge, green!50!black] (v0) -- node[midway, above, text=green!50!black, font=\footnotesize] {4} ++(-0.8, 0) node[left, text=black, font=\footnotesize] {$e_1$};
        \draw[myedge, red] (v0) to[out=120, in=0] node[midway, above, text=red, font=\footnotesize] {6} ++(-1.0, 0.5) node[left, text=black, font=\footnotesize] {$e_2$};
        \draw[myedge, red] (v1) to[out=120, in=0] node[midway, above, text=red, font=\footnotesize] {6} ++(-1.5, 0.8) node[left, text=black, font=\footnotesize] {$e_3$};

        \draw[myedge, green!50!black] (v8) -- node[midway, above, text=green!50!black, font=\footnotesize] {4} ++(0.8, 0) node[right, text=black, font=\footnotesize] {$e_1$};
        \draw[myedge, red] (v7) to[out=60, in=180] node[midway, above, text=red, font=\footnotesize] {6} ++(1.5, 0.8) node[right, text=black, font=\footnotesize] {$e_2$};
        \draw[myedge, red] (v8) to[out=60, in=180] node[midway, above, text=red, font=\footnotesize] {6} ++(1.0, 0.5) node[right, text=black, font=\footnotesize] {$e_3$};

    \end{tikzpicture}
    \par\smallskip
    {\small (c) Base strip of order 9.}
    \caption{Base and extension strips for the 3-ecd colouring of
    $C_n(1,2)$. Joining dangling edges with the same boundary name $e_i$
    closes either base; their numeric colour labels agree.
    Inserting panel (b) increases the order by 4.}
    \label{fig:circulant-12-construction}
\end{figure}

\clearpage

For a colour class $\mathcal C$, let its \emph{support} be the set of
vertices belonging to its cycles. The following observation isolates the
rotation argument used below.

\begin{lemma}[Merging by rotation]
\label{lem:circulant-rotation-merging}
Let a circulant $G$ of order $n$ be the edge-disjoint union of spanning
circulants $G_1,\ldots,G_r$. For each $i$, choose a $k_i$-ecd of $G_i$ and
select one of its colour classes; let $s_i$ be the size of its support. If
the graphs are indexed so that
\[
    n>s_i\sum_{j<i}s_j
    \qquad\text{for }i=2,\ldots,r,
\]
then
\[
    \ecd(G)\leq 1+\sum_{i=1}^r(k_i-1).
\]
\end{lemma}

\begin{proof}
We rotate the chosen ecds so that the selected colour classes have pairwise
disjoint supports; we can then merge these classes into one colour class.
For every $x\in\mathbb Z_n$, the map $v\mapsto v+x$ is an automorphism of
each $G_i$. Thus adding $x$ to every vertex of a chosen ecd merely rotates
its colour classes and does not affect the validity of the ecd. Suppose the
selected classes in
$G_1,\ldots,G_{i-1}$ already have pairwise disjoint supports, and let $U$ be
the union of these supports. Let $S_i$ be the support of the selected class
in $G_i$. For a fixed $u\in S_i$, exactly $|U|$ values of $x$ place $u+x$
in $U$. Hence at most
$|S_i||U|=s_i\sum_{j<i}s_j<n$ values of $x$ make the rotated support
intersect $U$. Hence at least one shift makes the two supports disjoint.
After all supports
have been made disjoint, merge the $r$ selected classes. Their union
is still a disjoint union of even cycles; every other colour class is left
unchanged.
\end{proof}

\begin{theorem}
\label{thm:circulant-1234}
For every odd $n\geq9$,
\[
    \ecd(C_n(1,2,3,4))=5.
\]
\end{theorem}

\begin{proof}
Proposition~\ref{prop:ecd-index-chromatic-index} gives
$\ecd(C_n(1,2,3,4))\geq5$.
For odd $n\geq13$, decompose the edge set into $C_n(1,3)$ and $C_n(2,4)$.
By Theorem~\ref{thm:circulant-12-13}, these factors have 3-ecds shown in
Figures~\ref{fig:circulant-13-construction} and
\ref{fig:circulant-12-construction}. In the first figure, the class formed
by colours 1 and 2 is supported within $\{0,1,2,3,4,5,6\}$, with support
$\{0,1,3,4\}$ when $n=13$; the extension strip adds no edges of this
class. In the second figure, the class formed by colours 1 and 2 is supported
on four consecutive vertices, and its extension strip adds no edges of this
class. Multiplication by 2 maps this latter class to $C_n(2,4)$ and its
support to a four-term progression of difference 2. For $n=13$, rotate this
progression to $\{5,7,9,11\}$; for odd $n\geq15$, rotate it to
$\{7,9,11,13\}$. These sets are disjoint from the corresponding supports
in $C_n(1,3)$. Merging the selected classes therefore gives a 5-ecd.

For $n=11$, the cycles $(0,1,5,2,3,10,9,6,7,8)$,
$(0,7,9,5,8,1,4,3,6,10)$, $(0,4,8,6,2,10,1,9)$,
$(0,2,4,5,7,3)$, and $(1,2,9,8,10,7,4,6,5,3)$ give a 5-ecd,
with one cycle in each class.
Each cycle closes from its last vertex to its first. Their edges partition
$E(C_{11}(1,2,3,4))$.
Finally, $C_9(1,2,3,4)=K_9$, whose ecd index is 5 by
Theorem~\ref{thm:complete-index}.
\end{proof}

Complete graphs are a special case of circulants; in particular,
$K_{4n+1}=C_{4n+1}(1,\ldots,2n)$.

\begin{theorem}
    For $n \geq 2$, the ecd index of $K_{4n+1}$ is $2n+1$.
    \label{thm:complete-index}
\end{theorem}
\begin{proof}
The solution of Alspach's conjecture \cite{bryant2014cycle} gives a
decomposition of $K_{4n+1}$ into $2n$ cycles of length $4n$ and one cycle
of length $2n$. Giving each cycle a separate colour yields a $(2n+1)$-ecd.
Proposition~\ref{prop:ecd-index-chromatic-index} gives the matching lower
bound $\ecd(K_{4n+1})\geq2n+1$.
\end{proof}

\begin{corollary}
\label{cor:paired-circulant-factors}
Let $r\geq1$, and let $C_n(j_1,\ldots,j_{2r})$ have odd order, and suppose its generators
can be paired so that every resulting 4-regular spanning circulant has a
3-ecd with one colour class supported on at most six vertices. If
$n>36(r-1)$, then
\[
    \ecd(C_n(j_1,\ldots,j_{2r}))=2r+1.
\]
\end{corollary}

\begin{proof}
Apply Lemma~\ref{lem:circulant-rotation-merging} with $s_i\leq6$ and
$k_i=3$ to obtain index at most $2r+1$. Proposition
\ref{prop:ecd-index-chromatic-index} gives the matching lower bound.
\end{proof}

For example, this pairing condition holds whenever every paired factor
satisfies the hypotheses of Theorem~\ref{thm:circulant-two-hamiltonian-cycles},
which provides the required colour class on four vertices.

Finally, $C_n(1,\ldots,t)$ is $2t$-vertex-connected
\cite{boesch1984circulants}. Thus Theorems~\ref{thm:circulant-12-13}
and~\ref{thm:circulant-1234} also give infinite maximally vertex-connected
families of 4- and 8-regular graphs with ecd index 3 and 5, respectively.

\section{Maximally vertex-connected families}
\label{sec:maximally-connected}

\begin{figure}[ht]
    \centering
    \begin{tikzpicture}[scale=1]
        \tikzset{
            vertex/.style={fill=black, circle, inner sep=1.5pt},
        }

        \node at (-1.5, 3) {\textbf{$G$}};
        \node[vertex, label=above:$v$] (v) at (-1.5, 1.5) {};

        \node[vertex, label=below:$c_1$] (c1_G) at (-3.5, -1) {};
        \node[vertex, label=below:$c_2$] (c2_G) at (-2.5, -1) {};
        \node (dots_c_G) at (-1.5, -1) {\dots};
        \node[vertex, label=below:$c_{d-1}$] (cdm1_G) at (-0.5, -1) {};
        \node[vertex, label=below:$c_d$] (cd_G) at (0.5, -1) {};

        \draw (v) -- (c1_G);
        \draw (v) -- (c2_G);
        \draw (v) -- (cdm1_G);
        \draw (v) -- (cd_G);

        \draw[->, thick] (1.1, 1) -- (2.1, 1);

        \node at (5, 3) {\textbf{$\expand(G)$}};

        \node[vertex, label=above:$u_1$] (u1) at (3.5, 2) {};
        \node[vertex, label=above:$u_2$] (u2) at (4.5, 2) {};
        \node (dots_u) at (5.5, 2) {\dots};
        \node[vertex, label=above:$u_{d-1}$] (ud) at (6.5, 2) {};

        \node[vertex, label=right:$v_1$] (v1) at (3.3, 0.5) {};
        \node[vertex, label=right:$v_2$] (v2) at (4.3, 0.5) {};
        \node (dots_v) at (5.3, 0.5) {\dots};
        \node[vertex, label=right:$v_{d-1}$] (vdm1) at (6.3, 0.5) {};
        \node[vertex, label=right:$v_d$] (vd) at (7.3, 0.5) {};

        \node[vertex, label=below:$c_1$] (c1_f) at (3.3, -1) {};
        \node[vertex, label=below:$c_2$] (c2_f) at (4.3, -1) {};
        \node (dots_c_f) at (5.3, -1) {\dots};
        \node[vertex, label=below:$c_{d-1}$] (cdm1_f) at (6.3, -1) {};
        \node[vertex, label=below:$c_d$] (cd_f) at (7.3, -1) {};

        \draw (u1) -- (v1);
        \draw (u1) -- (v2);
        \draw (u1) -- (vdm1);
        \draw (u1) -- (vd);

        \draw (u2) -- (v1);
        \draw (u2) -- (v2);
        \draw (u2) -- (vdm1);
        \draw (u2) -- (vd);

        \draw (ud) -- (v1);
        \draw (ud) -- (v2);
        \draw (ud) -- (vdm1);
        \draw (ud) -- (vd);

        \draw (v1) -- (c1_f);
        \draw (v2) -- (c2_f);
        \draw (vdm1) -- (cdm1_f);
        \draw (vd) -- (cd_f);

    \end{tikzpicture}
    \caption[Bipartite construction]{The construction of $\expand(G)$ from $G$ by
    replacing a vertex by a bipartite gadget.}
    \label{fig:fg-construction}
\end{figure}
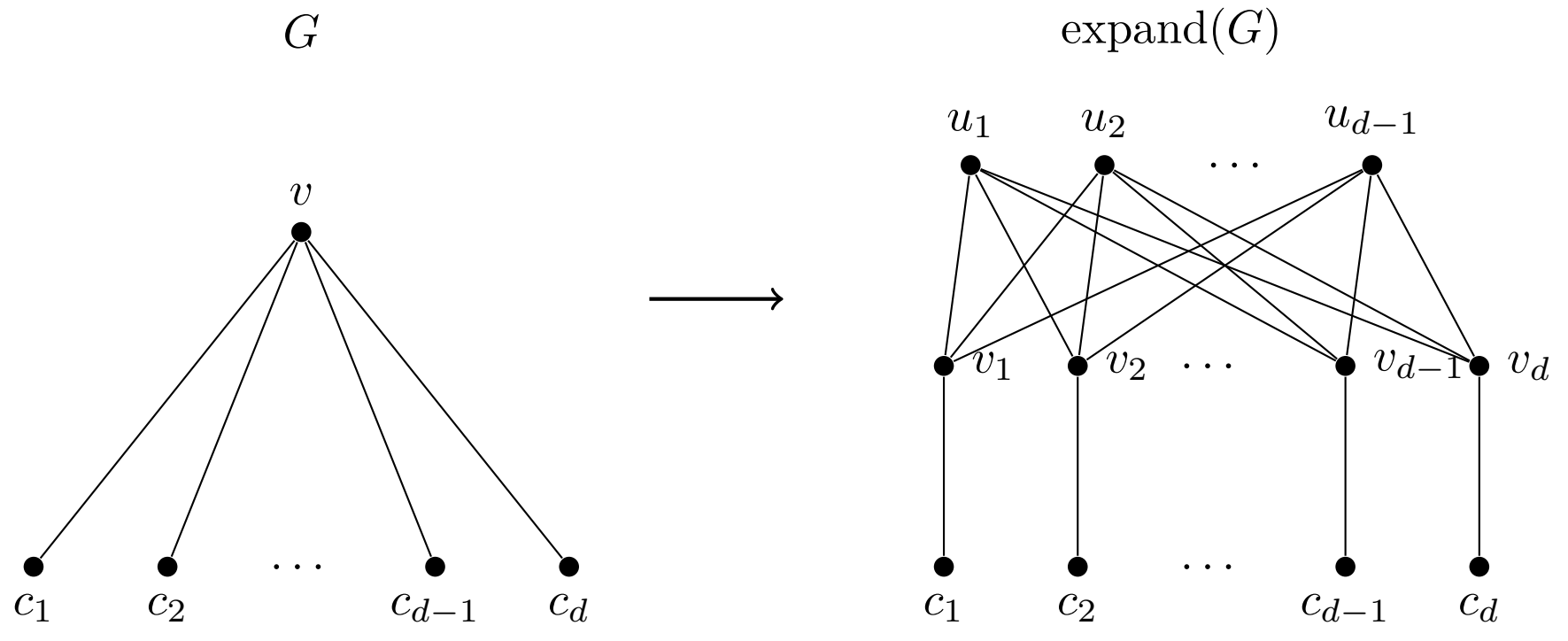

We next use Meredith's expansion operation~\cite{meredith1973regular} to obtain
maximally vertex-connected even-order families in every even degree at least 4.
Let $\expand$ denote an operation
that takes a finite $d$-regular ($d\geq1$)
graph $G$. For each $v \in V(G)$, list its incident edges as
$vc_1,\ldots,vc_d$, where the vertices $c_i$ need not be distinct.
Replace $v$ by $d$ vertices $v_1,\ldots,v_d$ and reattach $vc_i$ to $v_i$
for every $i$.
Moreover, we add $d-1$ new vertices $u_1, \ldots, u_{d-1}$.
Finally, we construct a complete bipartite graph between $v_1,
\ldots, v_{d}$ and $u_1, \ldots, u_{d-1}$.
If $G$ is a $d$-regular graph on $n$ vertices,
then $\expand(G)$ is a $d$-regular simple graph
on $n(2d-1)$ vertices. A visualization
of the construction is shown in Figure \ref{fig:fg-construction}.
The key property of this construction is that it transforms a multigraph into a simple graph while preserving edge-colouring properties.

\begin{proposition}[Edge-colourability under expansion]\textup{\cite[Theorem~2]{meredith1973regular}}
    Let $G$ be a $d$-regular
    graph. Then $\expand(G)$ is
    $d$-edge-colourable if and only if $G$ is $d$-edge-colourable.
    \label{prop:fg-edge-colouring}
\end{proposition}

\begin{proposition}[Connectivity under expansion]\textup{\cite[Theorem~3]{meredith1973regular}}
    Let $G$ be a finite $d$-regular $d$-edge-connected
    graph. Then $\expand(G)$ is $d$-vertex-connected.
    \label{prop:fg-vertex-connectivity}
\end{proposition}

\begin{theorem}[Expansion ecd bound]
    If $G$ is a connected $2k$-regular graph with an even number of edges,
    where $k\geq1$,
then $\expand(G)$ admits a $(k+1)$-ecd.
\label{thm:eulerian-expand-ecd}
\end{theorem}
\begin{proof}
Choose an Euler tour of $G$ and regard it as cyclic. At every vertex $v$,
the tour pairs the $2k$ incident edges into $k$ consecutive entrance--exit
pairs; at the chosen starting vertex, the last and first edges form one such
pair.

Label these pairs at $v$ by $1,\ldots,k$. If the two edges in the $i$-th
pair are reattached to $v_j$ and $v_\ell$, route the passage through the
gadget along the path $v_j u_i v_\ell$.
Every reattached edge is used once, every $v_j$ corresponds to a distinct edge
incident with $v$, and every $u_i$ is used for only one passage through $v$.
Thus the rerouted tour is a simple cycle $C$ of length $3\lvert E(G)\rvert$,
which is even by hypothesis. After removing $C$, in each gadget every $v_j$
and each of $u_1,\ldots,u_k$ has degree $2k-2$, while each of
$u_{k+1},\ldots,u_{2k-1}$ has degree $2k$. Thus the unused edges form a
bipartite graph $B$ of maximum degree at most $2k$, with all degrees even.
Lemma~\ref{lem:bipartite-ecd-index} gives it a $k$-ecd. Giving $C$ a new
colour yields a $(k+1)$-ecd of $\expand(G)$.
\end{proof}

\begin{theorem}[Exact ecd index after expansion]
Let $k\geq1$, and let $G$ be a connected $2k$-regular multigraph.
Then
\[
\ecd(\expand(G))=
\begin{cases}
\infty, & |E(G)| \text{ is odd},\\
k, & |E(G)| \text{ is even and $G$ is class 1},\\
k+1, & |E(G)| \text{ is even and $G$ is class 2}.
\end{cases}
\]
\label{thm:exact-expand-ecd}
\end{theorem}
\begin{proof}
The expansion has $(4k-1)|E(G)|$ edges, so an odd value of $|E(G)|$
makes an ecd impossible. Suppose that $|E(G)|$ is even. Theorem
\ref{thm:eulerian-expand-ecd} gives a $(k+1)$-ecd of the expansion.
By Proposition~\ref{prop:fg-edge-colouring}, the expansion is class 1 if
and only if $G$ is class 1. In that case
Proposition~\ref{prop:ecd-index-chromatic-index} gives ecd index $k$.
In the class 2 case, the same proposition rules out a $k$-ecd, so the
degree lower bound and the $(k+1)$-ecd give ecd index $k+1$.
\end{proof}

Expansion can decrease the ecd index. The $t=2$ multigraph of
Lemma~\ref{thm:three-vertex-multigraphs} has ecd index~6. Its odd order
rules out an 8-edge-colouring, and it has twelve edges, so
Theorem~\ref{thm:exact-expand-ecd} gives ecd index~5 after expansion.

\begin{lemma}
\label{lem:iterated-expansion}
Let $k\geq2$, and let $G$ be an even-order $2k$-regular,
$2k$-edge-connected multigraph that is not $2k$-edge-colourable. Repeatedly
applying $\expand$ to $G$ produces an infinite family of simple
$2k$-regular, $2k$-vertex-connected graphs of even order and ecd index
$k+1$.
\end{lemma}

\begin{proof}
Let $G_0=G$ and $G_{i+1}=\expand(G_i)$. By
Propositions~\ref{prop:fg-edge-colouring} and~\ref{prop:fg-vertex-connectivity},
$G_1$ is simple, $2k$-regular, $2k$-vertex-connected, and not
$2k$-edge-colourable. A $2k$-vertex-connected graph is also
$2k$-edge-connected, so the same argument applies at every iteration.
At each step the order is multiplied by $4k-1$, so it remains even and
strictly increases. Since $|E(G_i)|=k|V(G_i)|$ is even,
Theorem~\ref{thm:exact-expand-ecd} gives ecd index $k+1$ at every step.
\end{proof}

\begin{theorem}
\label{thm:maximally-vertex-connected-families}
For every $k\geq2$, there exists an infinite family of simple
$2k$-regular, $2k$-vertex-connected graphs of even order and ecd index
$k+1$.
\end{theorem}
\begin{proof}
For $k\neq3$, Meredith gives an even-order $2k$-regular,
$2k$-vertex-connected graph that is not $2k$-edge-colourable
\cite[pp.~57--59, Theorems~4 and~6 and the paragraph following
Theorem~6]{meredith1973regular}. For $k=3$, take the seed from Diko's
thesis~\cite[Chapter~4, Section~2]{diko2026thesis}.
Lemma~\ref{lem:iterated-expansion} applies in both cases.
\end{proof}

Thus even simple maximally vertex-connected regular graphs of even order
need not have the smallest possible ecd index $\Delta(G)/2$.
\section{Line graphs}
\label{sec:line-graphs}

Throughout this section, graphs are finite and the base graph $G$ is
simple. The line graph $L(G)$ is the graph
whose vertices correspond to the edges of $G$,
and two vertices of $L(G)$ are adjacent if and only if
their corresponding edges are adjacent in $G$.
The line graph of a $d$-regular graph is a $(2d-2)$-regular graph,
since each endpoint of an edge in $G$ is incident to $d-1$ other edges.
The \emph{oddness}
of a bridgeless cubic graph is the minimum number of odd cycles in a
2-factor \cite[Definition~3.9]{MARKSTROM20122676}.

Markstr\"om conjectured that the line graph of every 2-vertex-connected
cubic graph has an ecd \cite[Conjecture~3.1]{MARKSTROM20122676}.
He proved it for oddness~2 \cite[Theorem~3.10]{MARKSTROM20122676}, and
Liu, You, and Cui extended it to oddness~4
\cite[Theorem~4.7]{liu2020linegraphs}. The 2-vertex-connectivity is
necessary: if a cubic graph has a bridge, the corresponding vertex is a cut
vertex of its 4-regular line graph, and thus by
Proposition~\ref{prop:4reg-cut-vertex-no-ecd}, the line graph has no ecd.
As proved by Markstr\"om, the line graph of every 3-edge-colourable cubic
graph has an ecd \cite[Theorem~3.4]{MARKSTROM20122676}; Bonisoli and
Bonvicini later sharpened this to $\ecd(L(G))\leq 3$ for every class 1
cubic graph $G$ \cite[Corollary~9]{bonisoli2017even}.
They also sketch an extension of
their Proposition 8: for a class 1 $(2t+1)$-regular base graph with $t\geq1$, the
construction uses $3t$ subgraphs and gives a $3t$-ecd of its $4t$-regular
line graph \cite[Section~5, final two paragraphs]{bonisoli2017even}. The
bound $2t+1$ in Lemma~\ref{lem:line-graphs} is smaller when $t\geq2$ (and
equal when $t=1$), while Theorem~\ref{thm:exact-class1-line-graph-index}
gives the exact value $2t$ or $2t+1$, according to the parity of $|E(G)|$.
Diko \cite[Chapter~3, Theorem~3.7]{diko2026thesis} proved that $\ecd(L(G))\leq\binom{d}{2}$. Lemma~\ref{lem:line-graphs}
improves on this bound by grouping disjoint pairs through an edge-colouring of
$K_d$.

\begin{lemma}[Line-graph upper bound]
\label{lem:line-graphs}
Let $G$ be a simple $d$-regular graph with a proper $d$-edge-colouring, where
$d\geq2$. Then $\ecd(L(G))\leq d-1$ when $d$ is even and
$\ecd(L(G))\leq d$ when $d$ is odd.
\end{lemma}

\begin{proof}
Let $M_1,\ldots,M_d$ be the colour classes of a proper $d$-edge-colouring of
$G$. Since $G$ is $d$-regular, each $M_i$ is a perfect matching. For
$1\leq i<j\leq d$, the spanning subgraph $M_i\cup M_j$ is a union of even
cycles: at every vertex one incident edge has colour $i$ and one has colour
$j$, so the two colours alternate along every cycle. Let $F_{ij}$ be the
subgraph of $L(G)$ induced by the vertices corresponding to the edges in
$M_i\cup M_j$. It too is a union of even cycles. These cycles form an ecd
of $L(G)$ because each edge of $L(G)$ determines a unique colour pair.

It remains to group these cycles into colour classes. Regard the pair
$\{i,j\}$ as an edge of $K_d$. If $\{i,j\}$ and $\{a,b\}$ are disjoint,
then the edges of $G$ represented by the vertices of $F_{ij}$ have colours
$i,j$, whereas those represented by the vertices of $F_{ab}$ have colours
$a,b$. Hence $F_{ij}$ and $F_{ab}$ have disjoint vertex sets in $L(G)$,
even though the corresponding subgraphs of $G$ may meet at vertices of $G$.
Thus each colour class in a proper edge-colouring of $K_d$ gives one
colour class of the ecd. The bounds follow from $\chi'(K_d)=d-1$ for
even $d$ and $\chi'(K_d)=d$ for odd $d$.
\end{proof}

\begin{theorem}[Exact ecd index of class 1 line graphs]
\label{thm:exact-class1-line-graph-index}
Let $G$ be a connected simple $d$-regular graph with a proper $d$-edge-colouring, where
$d\geq2$. Then
\[
    \ecd(L(G))=
    \begin{cases}
        d-1,&\text{if $|E(G)|$ is even},\\
        d,&\text{if $|E(G)|$ is odd}.
    \end{cases}
\]
\end{theorem}

\begin{proof}
If $|E(G)|$ is even, Jaeger's theorem gives that $L(G)$ is class 1
\cite[Theorem~1.2.1]{jaeger1974indice}, so
Proposition~\ref{prop:ecd-index-chromatic-index} gives ecd index $d-1$.
If $|E(G)|$ is odd, then $d$ is odd because the $d$ colour classes are
perfect matchings of equal size. The $(2d-2)$-regular graph $L(G)$ then
has odd order, so the same proposition gives ecd index at least $d$.
Lemma~\ref{lem:line-graphs} gives the matching upper bound.
\end{proof}

For class 2 cubic base graphs, Markstr\"om's conjecture remains open.
Brinkmann et al.\ verified it computationally through order~36, as
reported by Liu, You, and Cui~\cite[Introduction]{liu2020linegraphs}.
Our order~38 computation found that every line graph of a cyclically
4-edge-connected snark of girth at least~5 has a $3$-ecd.

Among the class 2 4-regular graphs $G$ checked through order~16 and at
order~18, $L(G)$ has ecd index~4 exactly when $G$ has a cut vertex, and~3
otherwise~\cite[Chapter~3, subsection ``Line graphs'']{diko2026thesis}.
At a cut vertex $v$, the edges of a component of $G-v$ and its two edges
to $v$ give an odd part of a 4-edge cut in $L(G)$. This rules out a
$3$-ecd by Proposition~\ref{prop:edge-cut-no-d-edge-colouring}, and the
cut-vertex counts in Table~\ref{tab:ecd-4-regular}(b) account for all
computed cases of ecd index~4.

We also computed the ecd index of the line graphs of small 5-regular graphs.
For base graphs of orders 8 and 12, all the line graphs are class 1. For
orders 6, 10, and 14, all have ecd index~5~\cite[Chapter~3, subsection ``Line graphs'']{diko2026thesis}.
Recall that an \emph{$r$-graph} is an $r$-regular graph $G$ such that
$|\delta(S)|\geq r$ for every set $S\subseteq V(G)$ of odd cardinality.
Our computations suggest the following general form of Markstr\"om's
conjecture.
\begin{conjecture}
    For every $r\geq 3$, the line graph of every simple $r$-graph
    has an ecd.
\end{conjecture}

For $r=3$, the connected $r$-graphs are precisely the bridgeless cubic
graphs, so this includes Markstr\"om's conjecture. By
Lemma~\ref{lem:line-graphs}, it suffices to consider class 2 $r$-graphs.
The $r$-graph assumption cannot be
replaced by mere 2-vertex-connectivity or 2-edge-connectivity of the base
graph.
For every $k\geq1$, the complete graph $K_{4k+3}$ has vertex- and
edge-connectivity $4k+2$, yet its line graph has no ecd.
Indeed, $L(K_{4k+3})$ has $(4k+3)(2k+1)(4k+1)$ edges, an odd number.
The odd order of $K_{4k+3}$ also violates the $r$-graph condition for
$S=V(K_{4k+3})$.

\section*{Acknowledgements}

This work was partially supported from the research grants APVV-23-0076, VEGA 1/0173/25, and
VEGA 1/0613/26.

The authors used AI models from OpenAI, Anthropic, and Google to check and
improve proofs, edit text, and write portions of the code. The authors
reviewed and edited the resulting material and take full responsibility for
the article.

\bibliographystyle{abbrvurl}
\bibliography{literature}

\end{document}